\documentclass[a4paper,11pt]{article}
\usepackage{ifstyle}
\usepackage{graphicx}
\usepackage{epsfig} 
\numberwithin{equation}{section} 

\usepackage[english]{babel}               
\usepackage[T1]{fontenc}
\usepackage[latin1]{inputenc}
\usepackage{amssymb}   
\usepackage{array}                 
\def\Xint#1{\mathchoice
{\XXint\displaystyle\textstyle{#1}}%
{\XXint\textstyle\scriptstyle{#1}}%
{\XXint\scriptstyle\scriptscriptstyle{#1}}%
{\XXint\scriptscriptstyle\scriptscriptstyle{#1}}%
\!\int}
\def\XXint#1#2#3{{\setbox0=\hbox{$#1{#2#3}{\int}$}
\vcenter{\hbox{$#2#3$}}\kern-.5\wd0}}

\def\dashint{\Xint-}
\newtheorem{theorem}{Theorem} 
\newtheorem{corol}{Corollary}
\newtheorem{defini}{Definition}  
  
\newtheorem{lem}{Lemma}
\newtheorem{Claim}{Claim}

\def\opn#1#2{\def#1{\operatorname{#2} } } 
\opn\Rm{Rm}
\opn\Ric{Ric} 
\opn\Rc{Rc}
\opn\rk{rk}
\opn\Scal{Sc}
\opn\Tr{Tr}
\opn\Trac{Tr} 
\opn\det{det} 
\opn\diam{diam} 
\opn\dist{dist} 
\opn\div{div}
\opn\dim{dim}
\opn\Ker{Ker} 
\opn\loc{loc}
\opn\For{For}
\opn\Hil{Hil}
\opn\exp{exp}
\opn\Psh{Psh}
\opn\Vol{Vol}
\opn\exph{exph}
\opn\Herm{Herm}
\opn\End{End} 
\opn\Hess{Hess} 
\opn\Vol{Vol} 

\newcommand{\N}{\mathbb N}

\newcommand{\R}{\mathbb R}
\newcommand{\C}{\mathbb C}

\newcommand{\I}{\mathbb I}

\newcommand{\contract}{\mathrel{\kern-1.5pt\vrule width6.0pt height0.4pt depth0pt
                \vrule width0.4pt height4.0pt depth0pt}}
\newcommand{\retract}{\mathrel{\kern-1.5pt\vrule width0.4pt height4.0pt depth0pt
                          \vrule width6.0pt height0.4pt depth0pt}}

\newdimen\boxrulethickness \boxrulethickness=.07em
\newdimen\Openboxwidth \Openboxwidth=.72em
\newcommand{\Openbox}{%
 \leavevmode
 \hbox{%
   \hfil\vrule width\boxrulethickness
   \vbox to\Openboxwidth{%
     \advance\Openboxwidth -2\boxrulethickness
     \hrule height \boxrulethickness width\Openboxwidth\vfil
     \hrule height\boxrulethickness}%
   \vrule width\boxrulethickness\hfil
 }}

\begin{document}
\begin{center} 
\Large{\bf{On the $L^2$-extension of twisted holomorphic sections of singular hermitian line bundles.}}
\\
\vspace{0.4cm}
\large{Nefton Pali}
\end{center} 
\begin{abstract}
We prove a sharp Ohsawa-Takegoshi-Manivel type $L^2$-extension result for twisted holomorphic sections of singular hermitian line bundles over almost Stein manifolds. 
We establish as corollaries some extension results for pluri-twisted holomorphic sections of singular hermitian line bundles over projective manifolds.
\end{abstract}
\section{Introduction}
The $L^2$-extension problem that we will consider can be formulated in a vague way as follows. Let $Y$ be the zero set of a generically transverse holomorphic section $\sigma$ of a hermitian vector bundle $E$ over an almost Stein manifold and let $L$ be a hermitian line bundle satisfying suitable curvature conditions with respect to the curvature of $E$ and a quasi-plurisubharmonic weight $\varphi$. Then any holomorphic section of $K_X+L$ satisfying a certain $L^2$-condition with respect to $\varphi$ admits an extension with uniform  $L^2$-estimate.

The classic situation, which has been first considered in \cite{Oh-Ta}, and successively in \cite{Siu2}, \cite{Ber}, \cite{Mc-Va}, is when the hermitian bundle $E$ is the flat bundle $\C_X$ and $X$ is a bounded pseudoconvex domain. The situations considered in  \cite{Siu2} and \cite{Mc-Va} are reduced by the authors to this particular case.
In all this works it is sufficient to prove first the result in the case the weight function $\varphi$ is smooth and then to extract a weak limit solution. This is always possible by the uniformity of the $L^2$-estimate satisfied by the holomorphic extensions.

In this particular case there is no trouble with the regularising process of $\varphi$ since the curvature assumptions reduces to the pseudoeffectivness of $L$. We give a simple proof of this classic case in the appendix. However this is not needed in our proof of the general geometric case, so the reader not interested with this particular case can skip this part.

In the case $E$ is not hermitian trivial the lost of positivity in the curvature conditions occur even in the case $E$ is holomorphically trivial and $X$ is a ball. This is due to the presence of  the singular weight produced from the section $\sigma$ in the curvature hypothesis. In the case $r:=\rk_{_{\C}}E>1$ there is no way known so far, of regularising a quasi-plurisubharmonic function $\varphi$ by producing an arbitrary small lost of positivity of the current 
$$
T:=i\partial\bar\partial(\varphi+r\log|\sigma|^2)\,,
$$
even in restriction of the complete K\"{a}hler manifold $X\smallsetminus Y$. 
In the case $r=1$ the regularisation process produces an arbitrary small lost of positivity of the current $T$ all over $X\smallsetminus Y$. This is possible by the Lelong-Poincar\'{e} formula.


On the other hand by replacing $|\sigma|^2$ with $|\sigma|^2+\varepsilon^2$ in the definition of $T$ one loses a fixed amount of positivity over an open set $\{|\sigma|<\lambda\varepsilon\}$, depending on the curvature of $E$. This lost of positivity can be handled by applying a new perturbation method (inspired from the perturbation method  introduced in \cite{Dem1}) to the standard  $L^2$-theory.

This local picture is very close to the case when $X$ is a projective variety since by subtracting a hyperplane section $H$ to $X$ one can cover $X\smallsetminus H$ by an increasing family of coordinate pseudoconvex open sets. In the subsection \ref{Loc-Proj} we explain how our general positivity estimate simplifies in the particular local or projective case. Our proof of the general global case (see the main theorem \ref{Ext-tm} in section \ref{sng-ext-rst}) is based in large part on the previous works \cite{Man}, \cite{Dem3} in which the authors consider the global geometric context without the presence of the singularities. In our global and singular case there is an extra lost of positivity due to the singularities of the weight $\varphi$. By replacing the current $T$ with the current 
$$
i\partial\bar\partial\left(\varphi_{\varepsilon}+r\log(|\sigma|^2+\varepsilon^2)\right)\,,
$$
$\varphi_{\varepsilon}:=\log(e^{\varphi}+\varepsilon^2)$ we introduce a fixed lost of positivity over an arbitrary small neighborhood of the poles 
$$
\varphi^{-1}(-\infty)\cup \{|\sigma|=0\}\,.
$$
Our new perturbed $L^2$-method requires only the regularisation of the bounded quasi-plurisubharmonic functions $\varphi_{\varepsilon}$, which is a substantially simpler task than the general case considered in \cite{Dem1} and \cite{Dem2}.

An other difficulty
is that we can not let $\varepsilon\rightarrow 0$ directly in the process of weak limits extraction. The reason is due to the fact that the starting holomorphic extension produced by standard $L^2$-methods is not $L^2$ with respect to the weight $\varphi$. For this reason we need first to construct a holomorphic $L^2(e^{-\varphi})$-extension. This is possible since the almost Stein condition provides an arbitrary large amount of positivity $\tau\omega>0$ required to absorb the error therm produced by the extension obtained by local gluing. However the holomorphic $L^2(e^{-\varphi})$-extension $F'$ constructed in this way does not satisfy an uniform $L^2$-estimate. The estimate obtained in this way blows up when the shape of the domains involved in the exhaustion increases.

The extension with uniform $L^2$-estimate is constructed by applying the previous $\varepsilon$-weak limit extraction process with initial extension $F'$.

Allowing singularities requires a condition on the weight $\varphi$ with respect to $Y$. This condition is always satisfied for quasi-plurisubharmonic weights of analytic type plus a continuous rest. We wish to point out that the regularising process in \cite{Dem2} produces only approximations of analytic type plus a bounded rest, which can be made smooth only after blow-ups. Moreover in the case $\rk_{_{\C}}E>1$ the subtle game of positivity lost and weak limit extraction does not allow to approximate a general weight $\varphi$ in order to avoid the singularity condition with respect to $Y$.
As previously suggested the extension result in the case $X$ projective and $\rk_{_{\C}}E=1$ can be obtained (see subsection \ref{Loc-Proj}) by slightly modifying the arguments in \cite{Man}, \cite{Dem3}. In this and other particular situations (see subsection \ref{Loc-Proj}) the singularity condition on the weight is not needed.

In the last section we establish as corollaries some extension results for pluri-twisted holomorphic sections of singular hermitian line bundles over projective manifolds. 
We are able to perform our extensions from singular varieties under an integrability assumption. The technique we use has been invented by Siu \cite{Siu1}, \cite{Siu2} and drastically simplified in \cite{Pa1}.
\section{Basics of perturbed $L^2$-theory }
\subsection{The abstract existence result} 
We start by proving the following abstract existence theorem on Hilbert spaces (see also \cite{Hor}, \cite{Dem1}), which part (B) will be quite crucial for the rest of the paper.
\begin{theorem}\label{Abs-Ex-Hil}
Let $T:H_1\rightarrow H_2$, $S:H_2\rightarrow H_3$ be closed densely defined linear operators between Hilbert spaces such that 
$S\circ T=0$ and let $L:H_2\rightarrow H_2$ be a linear operator whith domain $D(L)=H_2$ such that for all $v\in D(T^*)\cap D(S)$ hold
\begin{eqnarray}\label{Abs-BKN}
 (Lv,v)\leq  \|T^*v\|^2+\|Sv\|^2\,.
\end{eqnarray}
{\bf (A)} Let $g\in \Ker S$ admitting a constant $C_g>0$ such that for all $v\in H_2$ hold
\begin{eqnarray}\label{Ex-Est}
|(g,v)|^2\leq C_g(Lv,v)\,.
\end{eqnarray}
Then there exist $u\in D(T)$ such that $Tu=g$ and $\|u\|^2\leq C_g$.
\\
\\
{\bf (B)} 
Let $Q:H_2\rightarrow H_2$ be a self-adjoint operator, $D(Q)=H_2$ such that $Q^2=Q$ and denote by $P:H_2\rightarrow \Ker S$ the orthogonal projection operator.
Let also  $g\in \Ker S$ and $C>0$, $\delta\geq 0$ constants such that  for all $v\in H_2$ hold
\begin{eqnarray}\label{Pert-Ex-Est}
|(g,v)|^2\leq C\left((L+\delta^2 Q)v,v\right)\,.
\end{eqnarray}
Then there exist $(u,h)\in D(T)\times H_2$ with norm $\|u\|^2+\|h\|^2\leq C$ such that $Tu+\delta PQ h=g$.
\end{theorem}
{\bf Remark.} In the case the linear operator $L$ is nonnegative, it admits an inverse and $g\in \Ker S\cap D(L^{-1})$ then one can choose $C_g=(L^{-1}g,g)$ in statement (A). In fact set $\langle \cdot,\cdot\rangle:=(L\,\cdot,\cdot)$. Then the Cauchy-Schwarz inequality implies
\begin{eqnarray*}
|(g,v)|^2=|\langle L^{-1}g,v\rangle |^2&\le&\langle L^{-1}g, L^{-1}g\rangle\cdot\langle v,v\rangle 
\\
\\
&=&(g,L^{-1}g)\cdot (Lv,v)=C_g(Lv,v)\,.
\end{eqnarray*}
The key point about statement (B) is that in the applications a few quantities depend on the parameters  $Q,\delta$ with the exception of $C>0$ which is always uniform in $Q,\delta$. This will allow to take weak limits in order to obtain a solution of the equation $Tu=g$ even when the estimate \eqref{Ex-Est} is not available.\hfill$\Box$
\\
\\
{\bf Proof of (A).} The fact that $S\circ T=0$ implies $(\Ker S)^{\perp}\subset \Ker T^*$, thus
$D(T^*)=(\Ker S\cap D(T^*))\oplus (\Ker S)^{\perp}$. The equation $Tu=g$ is equivalent to the condition $(v,Tu)=(v,g)$ for all $v\in  D(T^*)$, which in its turn is equivalent with
\begin{eqnarray}\label{L2ExPf}
(T^*v,u)=(v,g)\,, 
\end{eqnarray}
for all $v\in  D(T^*)$. By the other hand combining the conditions \eqref{Abs-BKN} and  \eqref{Ex-Est} we deduce the estimate
$$
|(g,v)|\leq C^{1/2}_g \|T^*v\|\,,
$$ 
for all $v\in \Ker S\cap D(T^*)$, thus for all $v\in  D(T^*)$. This implies that the map $F:R(T^*)\rightarrow \C$, $F(T^*v):=(v,g)$
is well defined on the range $R(T^*)$ of $T^*$, linear and bounded. Therefore there exist a unique $u\in \overline{R(T^*)}$ such that $F(T^*v)=(T^*v,u)$ for all $v\in D(T^*)$, which is precisely the condition \eqref{L2ExPf}. Moreover $\|u\|=\|F\|\leq C^{1/2}_g$.
\\
\\
{\bf Proof of (B).} We consider the closed and densely defined linear operator 
$$
\tilde T: H_1\oplus H_2\longrightarrow H_2\,,\qquad \tilde T(u,h)=Tu+\delta PQ h\,.
$$
Clearly $S\circ \tilde T=0$. We remark that the adjoint operator $\tilde T^*$ is given by $\tilde T^*(v)\equiv (T^*_1,T^*_2)(v)=(T^*v,\delta QP v)$, with $D(\tilde T^*)=D(T^*)$. In fact for all $u\in D(T)$, $v\in D(T^*)$ and $h\in H_2$ hold
\begin{eqnarray*}
(u,\tilde T^*_1 v)+(h,\tilde T^*_2 v)
&=&\left<(u,h),\tilde T^*v\right>
=\left<\tilde T(u,h),v\right>
\\
\\
&=&(u,T^* v)+(h,\delta\, QPv)\,.
\end{eqnarray*}
By combining the assumption on $Q$ with the estimate \eqref{Abs-BKN} we deduce for all $v\in  \Ker S \cap D(T^*)$ the estimate 
$$
\left((L+\delta^2 Q)v,v\right)\leq  \|T^*v\|^2+\|\delta \,Qv\|^2=\|\tilde T^*v\|^2\,.
$$
Then the conclusion follows from the proof of statement (A) by replacing $T$ with $\tilde T$. 
\hfill $\Box$
\\
\\
The following elementary claims will be very useful in the process of extraction of weak limits.
\begin{Claim}\label{Prod-w-lim}
Let $(u_l)_l$ be a uniformly bounded sequence of functions point-wise converging  to $u$ almost everywhere and $(F_l)_l\subset L^2(X)$ be a sequence of functions $L^2$-weakly convergent to $F$. Then the sequence $(u_lF_l)_l$  converges $L^2$-weakly to $uF$.
\end{Claim}
$Proof$. We write
$$
u_l\,F_l-u\,F=u_l\,F_l-u\,F_l+u\,F_l-u\,F=(u_l-u)F_l+u(F_l-F)\,.
$$
Obviously the sequence $u(F_l-F)$ converges to zero $L^2$-weakly. For any $\psi\in L^2(X)$ we remark the inequality
$$
\left|\int_X(u_l-u)F_l\,\psi\,dV\right |^2\leq \|F_l\|^2_{L^2(X)}\int_X|u_l-u|^2|\psi|^2\,dV\,.
$$
The $L^2$-weak convergence implies $\sup_l\|F_l\|^2_{L^2(X)}<+\infty$. The dominated convergence theorem implies $\int_X|u_l-u|^2|\psi|^2\,dV\rightarrow 0$. We deduce that also the sequence $(u_l-u)F_l$ converges to zero $L^2$-weakly.\hfill $\Box$
\begin{Claim}\label{w-lim-wei}
Let $(f_k)_k$, $(\varphi_k)_k$, $\varphi_k\leq c$ be two families of functions such that $\int_X|f_k|^2e^{-\varphi_k}\,dV\leq C$  and the sequence $\varphi_k$ point-wise converges to $\varphi$. 
\\
\\
{\bf (A).} 
Then there exist a $L^2(X,dV)$-weakly convergent subsequence $(f_l)_l$ with limit $f$, such that $\int_X|f|^2e^{-\varphi}\,dV\leq C$.
\\
\\
{\bf (B).} Assume also $\varphi_{k+1}\leq \varphi_k$, let $S\subset L^2(X,dV)$ be a closed subset, let 
$$
P_k:L^2(X,e^{-\varphi_k}dV)\longrightarrow S_k:=S\cap L^2(X,e^{-\varphi_k}dV)\,,
$$
be the $L^2(X,e^{-\varphi_k}dV)$-orthogonal projector and let $P$ the analogue projector on $S_{\infty}:=S\cap L^2(X,e^{-\varphi}dV)$ with respect to the weight $\varphi$. Then the subsequence $(P_lf_l)_l$ obtained by $(f_l)_l$ in $(A)$, converges $L^2(X,dV)$-weakly  to $Pf$.
\end{Claim}
$Proof\;of\;(A)$. By $L^2(X,dV)$-weak compactness there exist a $L^2(X,dV)$-weakly convergent subsequence  $F_l:=f_l\,e^{-\varphi_l/2}$ with weak limit $F$ such that $\int_X|F|^2\,dV\leq C$. By the other hand the uniformly bounded sequence $u_k:=e^{\varphi_k/2}$ point-wise converges to $u:=e^{\varphi/2}$. The conclusion follows from the fact that the sequence $f_l=u_l\,F_l$ converges  $L^2(X,dV)$-weakly to $f:=u\,F$ by claim \eqref{Prod-w-lim}. 
\\
$Proof\;of\;(B)$. The inequality $\int_X|P_lf_l|^2e^{-\varphi_l}\,dV\leq\int_X|f_l|^2e^{-\varphi_l}\,dV\leq C$ implies that one can extract a $L^2(X,dV)$-weak limit $g\in S_{\infty}$ from the sequence $(P_lf_l)_l$. Let $P^{\bot}_l$ be the complementary projectors. We infer $h_l:=P^{\bot}_lf_l\rightarrow h:=f-g$, $L^2(X,dV)$-weakly. We will be done if we prove the convergence
\begin{eqnarray}\label{cv-wlim}
0= \int_X\left<h_l,v\right>e^{-\varphi_l}\,dV\longrightarrow\int_X\left<h,v\right>e^{-\varphi}\,dV\,,
\end{eqnarray}
for all $v\in S_{\infty}\subset S_l$. The inequality $\int_X|h_l|^2e^{-\varphi_l}\,dV\leq C$ implies by the proof of statement (A) the  $L^2(X,dV)$-weak convergence $h_l\,e^{-\varphi_l/2}\rightarrow h\,e^{-\varphi/2}$. By the other hand the monotone convergence theorem implies the $L^2(X,dV)$-strong convergence
$v\,e^{-\varphi_l/2}\rightarrow v\,e^{-\varphi/2}$. We infer \eqref{cv-wlim}.
\hfill $\Box$
\subsection{The geometric context}
We explain now the geometric context related with the fundamental abstract inequality \eqref{Abs-BKN}. Let $(F,h)$ be a hermitian vector bundle over a Riemann manifold $(M,g)$ (of dimension $n$) equipped with a $h$-hermitian connection and let $\nabla$ be the induced hermitian connection over the hermitian vector bundle 
$((T^*_M)^{\otimes p}\otimes_{_\R}F, \left<\cdot,\cdot\right>)$, where $\left<\cdot,\cdot\right>$ is the induced hermitian product. We remind that the trace of any linear map $L:T_M\rightarrow T_M\otimes_{_\R}F$ is the element in $F$ identified by the map 
$$
\Tr L:\Lambda^nT_M\rightarrow \Lambda^nT_M\otimes_{_\R}F\,,\quad v_1\wedge...\wedge v_n\mapsto \sum_{j=1}^nv_1\wedge...\wedge L(v_j)\wedge...\wedge v_n\,.
$$
Then for any $F$-valued $2$-tensor $\alpha$, we define its trace $\Tr_g\alpha$ with respect to $g$ by 
$$
\Tr_g\alpha:=\Tr[(I_F\otimes g^{-1})\circ \alpha:T_M\rightarrow T_M\otimes_{_\R}F]\,.
$$
Moreover consider the first order differential operator
$$
\nabla:{\cal E}((T^*_M)^{\otimes p}\otimes_{_\R}F)\rightarrow 
{\cal E}((T^*_M)^{\otimes p+1}\otimes_{_\R}F),
$$
defined by $\nabla\alpha(\xi_0,...,\xi_p):=\nabla_{\xi_0}\alpha(\xi_1,...,\xi_p)$ for all vectors $\xi_0,...,\xi_p\in T_{M,x}$. The divergence of $\alpha$ is defined by the formula
$$
\div\alpha\,(\xi_1,...,\xi_{p-1}):=\Tr_g\nabla\alpha(\cdot,\cdot,\xi_1,...,\xi_{p-1}).
$$
The formal adjoint operator of $\nabla$ is given by the formula $\nabla^*=-\div$. In fact let $\beta$ be a $F$-valued $(p+1)$-tensor such that the intersection of its support whith the support of $\alpha$ is relatively compact and consider the $1$-form $\gamma$ given by the formula 
$\gamma(\xi)=\left<\alpha,\xi\contract \beta\right>$ for all $\xi\in T_M$. We prove first the identity
\begin{eqnarray}\label{div-forI}
\div\gamma=\left<\nabla\alpha,\beta\right>+\left<\alpha,\div\beta\right>.
\end{eqnarray}
In fact let $(e_k)_k$ be a local $g$-orthonormal frame of $T_M$ such that $\nabla {e_k}_{|_x}=0$ in a point $x$. Then at the point $x$ we have the equalities
\begin{eqnarray*}
\div\gamma_{|_x}&=&\sum_k\nabla_{e_k}\gamma (e_k)_{|_x}=\sum_ke_k\,.\gamma (e_k)_{|_x}=\sum_ke_k\,.\left<\alpha,e_k\contract \beta\right>_{|_x}
\\
&=&
\sum_k\left(\left<\nabla_{e_k}\alpha,e_k\contract \beta\right>_{|_x}+\left<\alpha,\nabla_{e_k}(e_k\contract \beta)\right>_{|_x}\right)
\\
&=&\left<\nabla\alpha,\beta\right>+\left<\alpha,\div\beta\right>_{|_x}.
\end{eqnarray*}
Let $\eta\in {\cal E}(M,T_M \otimes_{_\R}\C)$ be the unique vector field such that $\gamma=\eta\contract g$. By means of the partitions of unity we can assume that the support of $\gamma$ is contained in an oriented coordinate open set $U$. Let $V_g$ be the volume form over $U$, compatible with its orientation, induced by the metric $g$. With this notations hold the formula
\begin{eqnarray}\label{div-forII}
\div\gamma = \frac{d(\eta\contract V_g)}{V_g} \,.
\end{eqnarray}
(Remark that the expression on the right is independent of the orientation.)
By the Stokes formula we infer $\int_M (\div \gamma )\,dV_g=0$ which implies $\nabla^*=-\div$ by the identity \eqref{div-forI}. We prove now the formula \eqref{div-forII}. With the previous notations we consider the dual frame $e^*_k:=e_k\contract g$. Then $V_g=\pm e_1^*\wedge....\wedge e_n^*$. The fact that $\nabla {e^*_k}_{|_x}=0$ implies $d(e_k\contract V_g)_{|_x}=0$ by the covariant expression of the exterior differential (see formula \eqref{Ext-Dif-Cov} below with $F=\C\times M$). By differentiating the trivial identity
$$
\eta\contract V_g=\sum_k\gamma (e_k)\,(e_k\contract V_g)
$$
we infer at the point $x$
\begin{eqnarray*}
d(\eta\contract V_g)= \sum_{j,k}e_j\,.\gamma (e_k)\,e_j^*\wedge (e_k\contract V_g)=\sum_ke_k\,.\gamma (e_k)\,V_g=(\div \gamma )\,V_g\,.
\end{eqnarray*}
The $h$-hermitian connection $\nabla_F$ on $F$ extends to a exterior derivation on the sheaf ${\cal E}(\Lambda^pT^*_M\otimes_{_\R}F)$ that we still denote by $\nabla_F$. The relation with the operator $\nabla$ is
\begin{eqnarray}\label{Ext-Dif-Cov}
\nabla_F\,\alpha(\xi_0,...,\xi_p)=\sum_{j=0}^p(-1)^j\,\nabla\alpha(\xi_j,\xi_0,...,\hat\xi_j,...,\xi_p)\,.
\end{eqnarray}
In fact by expanding the right hand side therm we get
\begin{eqnarray*}
&&\sum_{0\leq j \leq p}(-1)^j\nabla_{F}(\alpha (\xi _0,...,\widehat{\xi _j},..., \xi _p))(\xi _j)
\\
\\
&-&\sum_{\scriptstyle 0\leq j,l \leq p
\atop
\scriptstyle j\not=l}
(-1)^j\alpha (\xi _0,...,\nabla_{\xi _l} \xi _j,...,\widehat{\xi _l},..., \xi _p)\,.
\end{eqnarray*}
The last sum is equal to the quantity
\begin{eqnarray}
&\displaystyle{
-\sum_{0\leq j<l \leq p}(-1)^{j+l}\alpha
(\nabla_{\xi _l}\xi _j  ,\xi _0,...,\widehat{\xi _j},...,\widehat{\xi _l},...,\xi _p)}\nonumber
\\\nonumber
\\
&\displaystyle{
+\sum_{0\leq l<j \leq p}(-1)^{j+l}\alpha
(\nabla_{\xi _l}\xi _j  ,\xi _0,...,\widehat{\xi _l},...,\widehat{\xi _j},...,\xi _p)}&\nonumber
\\\nonumber
\\
&\displaystyle{=\sum_{0\leq l<j \leq p}(-1)^{j+l}\alpha
(\nabla_{\xi _l}\xi _j -\nabla _{\xi _j}\xi _l  ,\xi _0,...,\widehat{\xi _l},...,\widehat{\xi _j},...,\xi _p)   }&\nonumber
\end{eqnarray}
The fact that the Levi-Civita connexion is torsion free allows to conclude.
We remark now that $\nabla_F^*=\nabla^*$ restricted to ${\cal E}(\Lambda^pT^*_M\otimes_{_\R}F)$. In fact this follows from the identity 
$\left<\nabla_F\,\alpha,\beta\right>=\left<\nabla\alpha,\beta\right>$, for all $F$-valued $p$-form $\alpha$ and $F$-valued $(p+1)$-form $\beta$. Let prove this identity. Let $(\theta_s)_s$ be a $h$-orthonormal frame of the bundle $F$ of complex rank $r$. The coefficients of the local expressions
$$
\nabla\alpha=\sum_{s=1}^r\sum_{j=1}^n\sum_{|I|=p}C_{j,I}^s\,e_j^*\otimes e^*_I\otimes \theta_s\,,\qquad \nabla_F\,\alpha=\sum_{s=1}^r\sum_{|K|=p+1}B_{K}^s\, e^*_K\otimes \theta_s\,,
$$
are related by the formula 
\begin{eqnarray*}
B_{K}^s=\sum_{j=0}^p(-1)^jC^s_{k_j, \hat K_j}\,,
\end{eqnarray*}
where $\hat K_j:=(k_0,...,\hat k_j,...,k_p)$. By the other hand
\begin{eqnarray*}
\left<\nabla\alpha,e^*_K\otimes \theta_s\right> &=&\sum_{j=1}^n\left<\nabla_{e_j}\alpha,e_j\contract e^*_K\otimes \theta_s\right>
=\sum_{j=1}^p\left<\nabla_{e_{k_j}}\alpha,e_{k_j}\contract e^*_K\otimes \theta_s\right>
\\
&=&
\sum_{j=1}^p(-1)^j\left<\nabla_{e_{k_j}}\alpha, e^*_{\hat K_j}\otimes \theta_s\right>=\sum_{j=0}^p(-1)^jC^s_{k_j, \hat K_j}=B_{K}^s
\\
&=&\left<\nabla_F\,\alpha,e^*_K\otimes \theta_s\right>\,,
\end{eqnarray*}
which proves the required identity. In conclusion we have found the formula 
\begin{eqnarray}\label{Adj-for}
\nabla_F^*\,\alpha=-\Tr_g\nabla\alpha\,,
\end{eqnarray}
for all $F$-valued $p$-forms $\alpha$.
We assume now that the hermitian vector bundle $(F,h)$ is over a K\"{a}hler manifold $(X,\omega)$ of complex dimension $n$. 
In local complex coordinates we will use the expression $\omega=\frac{i}{2}\sum_{k,l}\omega_{k,\bar l}\,dz_k\wedge d\bar z_l$.
Let $\partial_F$ and $\bar\partial_F$ be respectively the type $(1,0)$ and $(0,1)$ exterior derivatives induced by $\nabla_F$. Then the formal adjoint operator $\nabla^*_F$ splits as  $\nabla^*_F=\partial^*_F+\bar\partial^*_F$, where the formal adjoint $\partial^*_F$ of $\partial_F$ is a $(-1,0)$-degree operator and the formal adjoint $\bar\partial^*_F$ of $\bar\partial_F$ is a $(0,-1)$-degree operator. If $A$ and $B$ are two operators on the fibres of $\Lambda^{\bullet}T^*_M\otimes_{_\R}F$ then their bracket is defined by the formula
$$
[A,B]:=AB-(-1)^{\deg A\cdot \deg B}BA\,.
$$
The dual metric  
$$
\omega^*=2i\,\sum_{k,l}\omega^{l\bar{k}}\frac{\partial }{\partial z_k}\wedge\frac{\partial }{\partial \bar{z}_l}
$$
can be considered as an operator on the fibres of bundle $\Lambda^{\bullet}T^*_M\otimes_{_\R}F$. In fact we set $\omega^*\,\alpha:=\omega^*\contract \alpha$, where the contraction operator is defined by $(\xi\wedge \bar\eta)\contract \alpha:=\xi\contract(\bar\eta\contract \alpha)=-\alpha(\xi\wedge \bar\eta,\cdot)$, for all $\xi, \eta\in T_X^{1,0}$.
With this notations hold the basic K\"{a}hler identities
\begin{eqnarray}\label{Bas-Kah-id}
\partial^*_F=i[\omega^*,\bar\partial_F] \,,\qquad  \bar\partial^*_F=-i[\omega^*,\partial_F] \,.
\end{eqnarray}
In fact let $(z_1,...,z_n)$ be holomorphic $\omega$-geodesic coordinates centered at a point $x$ and set
$e_k:=\frac{\partial}{\partial x_k}{\vphantom{x} } _{| _{_x}}$, $\zeta_k:=e_k^{1,0}=\frac{\partial}{\partial z_k}{\vphantom{z} } _{| _{_x}}$ . Let
$$
\alpha=\sum_{\scriptstyle |K|=p
\atop
\scriptstyle |L|=q} \alpha_{K,L}\otimes dz_K\wedge d\bar z_L\,, 
$$ 
be a $(p,q)$-form and set $C_{K,L,r}:=\nabla_{F,\zeta_r}\alpha_{K,L}\in F_x$, $C_{K,L,\bar r}:=\nabla_{F,\bar\zeta_r}\alpha_{K,L}\in F_x$.
Then at the point $x$ hold the local expressions
\begin{eqnarray*}
\partial_F\alpha
& = &
\sum_{r=1}^n\sum_{\scriptstyle |K|=p
\atop
\scriptstyle |L|=q}C_{K,L,r}\otimes \zeta^*_r\wedge \zeta^*_K \wedge \bar\zeta^*_L\,,
\\
\\
\bar\partial_F\alpha
&=&
\sum_{r=1}^n\sum_{\scriptstyle |K|=p
\atop
\scriptstyle |L|=q}(-1^p)C_{K,L,\bar r}\otimes  \zeta^*_K \wedge \bar\zeta^*_r\wedge\bar\zeta^*_L\,.
\end{eqnarray*}
The fact that the family of real vectors $(e_k,Je_k)_k$ is a $g(x):=\omega(\cdot,J\cdot)(x)$-orthonormal base implies, by the formula \eqref{Adj-for}, the expansion at the point $x$
\begin{eqnarray*}
&&\nabla_F^*\,\alpha 
=
-\Tr_g\nabla\alpha\,=\,-\sum_{j=1}^n\left(e_j\contract \nabla _{e_j}\alpha+Je_j\contract \nabla _{Je_j}\alpha\right) 
\\
\\
&=&
-2\sum_{j=1}^n\left(\bar \zeta_j\contract \nabla _{\zeta_j}\alpha+\zeta_j\contract \nabla _{\bar\zeta_j}\alpha\right)
\\
\\
&=&
-2\sum_{j=1}^n\sum_{\scriptstyle |K|=p
\atop
\scriptstyle |L|=q}\left((-1^p)C_{K,L,j}\otimes \zeta^*_K\wedge (\bar \zeta_j\contract \bar\zeta^*_L)+C_{K,L,\bar j}\otimes (\zeta_j\contract \zeta^*_K)\wedge \bar\zeta^*_L\right)
\end{eqnarray*}
By bidegree reasons we infer at the point $x$ the expressions
\begin{eqnarray}
\partial^*_F\alpha
& = &
-2\sum_{j=1}^n\sum_{\scriptstyle |K|=p
\atop
\scriptstyle |L|=q}C_{K,L,\bar j}\otimes (\zeta_j\contract \zeta^*_K) \wedge \bar\zeta^*_L\,,\label{locAdj-par}
\\\nonumber
\\
\bar\partial^*_F\alpha
&=&
-2\sum_{j=1}^n\sum_{\scriptstyle |K|=p
\atop
\scriptstyle |L|=q}(-1^p)C_{K,L, j}\otimes  \zeta^*_K \wedge (\bar \zeta_j\contract\bar\zeta^*_L)\,.\label{locAdj-dpar}
\end{eqnarray}
Moreover
\begin{eqnarray*}
\omega^*\contract \alpha
=
2i(-1)^p\sum_{j=1}^n\sum_{\scriptstyle |K|=p
\atop
\scriptstyle |L|=q}\alpha_{K,L}\otimes \left(\frac{\partial}{\partial z_j}\contract dz_K\right)\wedge \left(\frac{\partial}{\partial \bar z_j}\contract d\bar z_L\right)\,+O(|z|^2)\,.
\end{eqnarray*}
We infer at the point $x$ the expressions
\begin{eqnarray*}
-i\bar\partial_F(\omega^*\contract \alpha)
&=&
-2\sum_{j,r=1}^n\sum_{\scriptstyle |K|=p
\atop
\scriptstyle |L|=q}C_{K,L, \bar r}\otimes (\zeta_j\contract \zeta^*_K) \wedge \bar\zeta^*_r\wedge (\bar \zeta_j\contract\bar\zeta^*_L)\,,
\\
\\
i\omega^*\contract\bar\partial_F\alpha
&=&
-2\sum_{j,r=1}^n\sum_{\scriptstyle |K|=p
\atop
\scriptstyle |L|=q}C_{K,L, \bar r}\otimes (\zeta_j\contract \zeta^*_K) \wedge (\bar \zeta_j\contract \bar\zeta^*_r\wedge\bar\zeta^*_L)\,.
\end{eqnarray*}
Then the identity 
\begin{eqnarray}\label{contrac-id}
\bar \zeta_j\contract \bar\zeta^*_r\wedge\bar\zeta^*_L=\delta_{j,r} \bar\zeta^*_L - \bar\zeta^*_r\wedge (\bar \zeta_j\contract\bar\zeta^*_L)\,,
\end{eqnarray}
combined with \eqref{locAdj-par} implies the first basic K\"{a}hler identity in \eqref{Bas-Kah-id}. By the other hand the expressions at the point $x$
\begin{eqnarray*}
i\partial_F(\omega^*\contract \alpha)
&=&
-2(-1)^p\sum_{j,r=1}^n\sum_{\scriptstyle |K|=p
\atop
\scriptstyle |L|=q}C_{K,L,  r}\otimes \zeta^*_r\wedge (\zeta_j\contract \zeta^*_K) \wedge  (\bar \zeta_j\contract\bar\zeta^*_L)\,,
\\
\\
-i\omega^*\contract\partial_F\alpha
&=&
-2(-1)^p\sum_{j,r=1}^n\sum_{\scriptstyle |K|=p
\atop
\scriptstyle |L|=q}C_{K,L, r}\otimes \left(\zeta_j\contract(\zeta^*_r \wedge \zeta^*_K)\right) \wedge (\bar \zeta_j\contract \bar\zeta^*_L)\,,
\end{eqnarray*}
combined with the conjugate of \eqref{contrac-id} and with \eqref{locAdj-dpar} implies the second basic K\"{a}hler identity in \eqref{Bas-Kah-id}.
\\
We set now by $(u,v )_{\omega,h}:=\int_X\left<u,v\right>_{\omega,h}\,dV_{\omega}$ the $L^2$-product and by  $\|v\|_{\omega,h}$ the corresponding norm. For readers convenience we give a proof of the following fundamental result \cite{Boc}, \cite{Kod}, \cite{Nak}, \cite{Dem1}.
\begin{theorem}{\bf($L^2$-Bochner-Kodaira-Nakano Inequality).} Let $(F,\bar\partial_F ,h)$ be a holomorphic hermitian vector bundle over a complete K\"{a}hler manifold $(X,\omega)$ such that $-C\omega\otimes \I_{F}\le i\,{\cal C}_h(F)$ for some constant $C>0$. Then the $L^2_{\omega,h}$-extension $($in the sense of distributions$)$ of the formal adjoint operator $\bar\partial_F^*$ coincides whith the Hilbert adjoint of the $L^2_{\omega,h}$-extension of $\bar\partial_F$ and
\begin{eqnarray}\label{Geom-BKN}
\left([i\,{\cal C}_h(F), \omega^*]v,v\right)_{\omega, h}\leq \|\bar\partial_F v\|^2_{\omega,h}+\|\bar\partial_F^*v\|^2_{\omega,h}\,.
\end{eqnarray}
for any $v\in D(\bar\partial_F)\cap D(\bar\partial_F^*)\subset  L^2_{\omega, h}(X,\Lambda^{p,q}T^*_X\otimes F)$.
\end{theorem}
$Proof$. The Chern connection of $F$ is defined by $\nabla_{F,h}:=\partial_{F,h}+\bar\partial_F$, with $\partial_{F,h}:=h^{-1}\cdot \partial_{\overline{F}^*}\cdot h$. We infer $\partial_F^2=0$. Thus the Chern curvature of the bundle $F$ satisfies the identity ${\cal C}_h(F)=[\partial_{F,h}\,,\, \bar\partial_F]$. Combining this with the basic K\"{a}hler identities \eqref{Bas-Kah-id} and with the Jacobi identity we obtain the Bochner-Kodaira-Nakano identity $\Delta''_F=\Delta'_F+[i\,{\cal C}_h(F), \omega^*]$. By integrating by parts we infer
\begin{eqnarray*}
\|\bar\partial_F v\|^2_{\omega,h}+\|\bar\partial_F^*v\|^2_{\omega,h}
&=&
(\Delta''_F v,v )_{\omega,h}\,=\,(\Delta'_Fv,v )_{\omega,h}+\left([i\,{\cal C}_h(F), \omega^*]v,v\right)_{\omega, h}
\\
\\
&=&
\|\partial_F v\|^2_{\omega,h}+\|\partial_F^*v\|^2_{\omega,h}+\left([i\,{\cal C}_h(F), \omega^*]v,v\right)_{\omega, h}
\\
\\
&\ge&
\left([i\,{\cal C}_h(F), \omega^*]v,v\right)_{\omega, h}\,.
\end{eqnarray*}
for all $v\in {\cal D}^{p,q}(X,F)$. We extend now this inequality to $L^2$-sections. For the moment we denote by $\bar\partial_{F, fm}^*$  the $L^2_{\omega,h}$-extension (in the sense of distributions) of the formal adjoint operator $\bar\partial_F^*$ and by 
$\bar\partial_{F, Hb}^*$ the Hilbert adjoint of the $L^2_{\omega,h}$-extension of $\bar\partial_F$. We remark first
\begin{eqnarray*}
D(\bar\partial_{F, fm}^*)&=&\left\{u\in L^2_{\omega,h}\,|\,\exists g\in L^2_{\omega,h}\,:\,(\varphi,g)=(\bar\partial_F\varphi,u)\,,\,\forall \varphi\in {\cal D}\right\} \,,
\\
\\
D(\bar\partial_{F, Hb}^*)&=&\left\{u\in L^2_{\omega,h}\,|\,\exists g\in L^2_{\omega,h}\,:\,(v,g)=(\bar\partial_Fv,u)\,,\,\forall v\in D(\bar\partial_F)\right\} \,.
\end{eqnarray*}
The obvious inclusion ${\cal D}\subset D(\bar\partial_F)$ implies $D(\bar\partial_{F, Hb}^*)\subset D(\bar\partial_{F, fm}^*)$. The identity $\bar\partial_{F, fm}^*=\bar\partial_{F, Hb}^*$ will follow immediately from the fact that for all $v\in D(\bar\partial_F)$ there exist a sequence $(\varphi_k)_k\subset {\cal D}$ such that $\varphi_k\rightarrow v$ and $\bar\partial_F\varphi_k\rightarrow \bar\partial_Fv$ in the $L^2_{\omega,h}$-norm. In fact the completeness assumption is equivalent to the existence of a non-decreasing sequence of functions $\chi_k\in {\cal D}(X,[0,1])$ such that the family of compact sets $\chi_k^{-1}(1)$ covers $X$ and $|d\chi_k|_{\omega}\leq 2^{-k}$. Thus $\chi_kv\rightarrow v$ and 
$$
\bar\partial_F(\chi_kv)=\chi_k\bar\partial_Fv+\bar\partial\chi_k\wedge v\longrightarrow \bar\partial_Fv\,,
$$ 
in the $L^2_{\omega,h}$-norm. This allows to localise the problem. The conclusion follows by using partitions of unity and usual regularising kernels $(\rho_{\varepsilon})$. From now on we identify the notations $\bar\partial^*_F\equiv \bar\partial_{F, fm}^*=\bar\partial_{F, Hb}^*$.
In the particular case $v\in D(\bar\partial_F)\cap D(\bar\partial^*_F)$ the sequence $(\varphi_k)_k\subset {\cal D}$ so far constructed satisfies also the extra condition $\bar\partial^*_F\varphi_k\rightarrow \bar\partial^*_Fv$ in the $L^2_{\omega,h}$-norm. This combined with the inequality
\begin{eqnarray*}
\left([i\,{\cal C}_h(F)+2C\omega, \omega^*]\varphi_k,\varphi_k\right)_{\omega, h}\leq \|\bar\partial_F \varphi_k\|^2_{\omega,h}+\|\bar\partial_F^*\varphi_k\|^2_{\omega,h}+2Cq\| \varphi_k\|^2_{\omega,h}\,,
\end{eqnarray*}
implies the required $L^2$-Bochner-Kodaira-Nakano Inequality. In order to prove this convergence we remark as before that
$$
\bar\partial^*_F(\chi_kv)=\chi_k\bar\partial^*_Fv-(\bar\partial\chi_k)^* v\longrightarrow \bar\partial^*_Fv\,,
$$ 
in the $L^2_{\omega,h}$-norm. So we can assume that the support of $v$ is contained in a coordinate open set trivialising the bundle $F$. The conclusion will follows from the fact that
\begin{eqnarray}\label{Fredc}
\|(\bar\partial_F^*v)\star \rho_{\varepsilon}-\bar\partial_F^*(v\star \rho_{\varepsilon}) \|^2_{\omega,h}\longrightarrow 0\,,
\end{eqnarray}
as $\varepsilon\rightarrow  0$. This is a straightforward consequence of the Friedrich's lemma.
\hfill$\Box$
\begin{lem}\label{O-T-Key-Es}
Let $\alpha\in \Lambda^{1,0}T^*_X\;,\;u\in \Lambda^{n,q}T^*_X\otimes_{_{\C}}F\;,\;v\in \Lambda^{n,q+1}T^*_X\otimes_{_{\C}}F$. Then	
\begin{eqnarray}
|\left<\bar \alpha\wedge u,v\right>_{\omega,h}|^2\le  |u|^2_{\omega,h}\left<[i\alpha\wedge \bar \alpha, \omega^*]v,v\right>_{\omega,h}\,.
\end{eqnarray}
\end{lem}
$Proof$. We prove first the identity $i\bar \alpha^*=[\alpha, \omega^*]$. Let $(\zeta_k)_k$ be a frame such that $\omega=i\sum_k\zeta^*_k\wedge \bar\zeta^*_k$. We set for notation simplicity $\zeta^*_{K,\bar L}:=\zeta^*_K\wedge \bar\zeta^*_L$ and write $\alpha=\sum_j\alpha_j\,\zeta^*_j$. The identity
\begin{eqnarray*}\label{exp-OT}
\left<\zeta^*_{K,\bar L}\,,\, \bar \alpha^*\zeta^*_{K,\bar H}\right>&=&\sum_j\,\left<\bar\alpha_j\,\bar\zeta^*_j\wedge \zeta^*_{K,\bar L}\,,\, \zeta^*_{K,\bar H}\right>\nonumber
\\
&=&(-1)^{|K|}\sum_j\,\left< \zeta^*_K\wedge \bar\zeta^*_j\wedge\bar\zeta^*_L\,,\, \alpha_j\,\zeta^*_{K,\bar H}\right>\,,
\end{eqnarray*}
combined with $\left<\bar\zeta^*_j\wedge \bar\zeta^*_L\,, \bar\zeta^*_H\right>=\left<\bar\zeta^*_L\,,\, \bar\zeta_j\contract \bar\zeta^*_H\right>$ implies
$$
i\bar \alpha^*\zeta^*_{K,\bar H}=(-1)^{|K|}i\sum_j\alpha_j\,\zeta^*_K\wedge (\bar\zeta_j\contract \bar\zeta^*_H)\,.
$$
We expand the therm
\begin{eqnarray*}
[\alpha, \omega^*]\zeta^*_{K,\bar H}&=&
\alpha\wedge(-1)^{|K|}i\sum_j\,(\zeta_j\contract \zeta^*_K)\wedge (\bar\zeta_j\contract \bar\zeta^*_H)
-
\omega^*\contract \sum_r\alpha_r\,\zeta^*_r\wedge\zeta^*_{K,\bar H}
\\
&=&
(-1)^{|K|}i\sum_{r,j}\alpha_r\,\zeta^*_r\wedge(\zeta_j\contract \zeta^*_K)\wedge (\bar\zeta_j\contract \bar\zeta^*_H)
\\
&+&
(-1)^{|K|}i\sum_{r,j}\alpha_r\,\zeta_j\contract (\zeta^*_r\wedge\zeta^*_K)\wedge (\bar\zeta_j\contract \bar\zeta^*_H)\;=\;i\bar \alpha^*\zeta^*_{K,\bar H}\,,
\end{eqnarray*}
by the conjugate contraction identity \eqref{contrac-id}.
By the other hand
$$
[i\alpha\wedge \bar \alpha, \omega^*]v=i\alpha\wedge \bar \alpha\wedge ( \omega^*\contract v)=-i\bar \alpha\wedge \alpha\wedge( \omega^*\contract v)=\bar \alpha\wedge\bar \alpha^*v\,,
$$
thus 
\begin{eqnarray}\label{nor-cont}
|\alpha^*v|^2=  \left<[i\alpha\wedge \bar \alpha, \omega^*]v,v\right>\,.
\end{eqnarray}
Then the conclusion follows by the Cauchy-Schwarz inequality
$$
|\left<\bar \alpha\wedge u,v\right>|^2=|\left< u,\bar \alpha^*v\right>|^2\le  |u|^2|\alpha^*v|^2\,.
$$
\hfill $\Box$
\\
We explain now the crucial inequality required for the proof of the extension results. This is due to \cite{Oh-Ta}, based on the previous works of \cite{Do-Fe}, \cite{Do-Xa}. The original proof in  \cite{Oh-Ta} has been substantially simplified in \cite{Oh}. The computation in this paper indicates that this inequality can also be obtained perturbing by a conformal factor the hermitian metric in the classic Bochner-Kodaira-Nakano Inequality. This has been also remarked recently in \cite{Pa2}.
\begin{lem}{\bf(Perturbed $L^2$-Bochner-Kodaira-Nakano Inequality).} \label{Prt-BKN}
\\
Let $(F,\bar\partial_F ,h)$ be a holomorphic hermitian line bundle over a complete K\"{a}hler manifold $(X,\omega)$, let $\varepsilon>0$ be a constant, let $\eta,\lambda>\varepsilon$ be two bounded smooth functions such that
$$
-C\omega\le \Theta^{\eta,\lambda}_h:=\eta\,i\,{\cal C}_h(F)-i\partial\bar\partial\eta-\lambda^{-1}i\,\partial\eta\wedge \bar\partial\eta\,,
$$
$
i\,\partial\eta\wedge \bar\partial\eta \le C\omega
$,
for some constant $C>0$
and set 
$
L_{\omega, h}^{\eta,\lambda}:=[\Theta^{\eta,\lambda}_h ,\, \omega^*]
$.
Then 
\begin{eqnarray}\label{Prt-Geom-BKN}
\left(L_{\omega, h}^{\eta,\lambda}\,v,v\right)_{\omega, h}\leq \|\bar\partial_F v\|^2_{\omega,\eta h}+\|\bar\partial_F^*v\|^2_{\omega,(\eta+\lambda)h}\,.
\end{eqnarray}
for any $v\in D(\bar\partial_F)\cap D(\bar\partial_F^*)\subset  L^2_{\omega, h}(X,\Lambda^{n,q}T^*_X\otimes F)$.
\end{lem}
$Proof$.
Let $\bar\partial_{F,\eta}^*$ be the Hilbert adjoint operator of $\bar\partial_F$ with respect to the conformal hermitian metric $\eta h$. By definition we infer the formula 
\begin{eqnarray*}
\bar\partial_{F,\eta}^*= \bar\partial_F^*-\eta^{-1}(\bar\partial \eta)^*\,.
\end{eqnarray*}
In particular $D(\bar\partial_{F,\eta}^*)=D(\bar\partial_F^*)$ by the existence of the constants $\varepsilon$ and $C$.
Moreover 
$$
i\,{\cal C}_{\eta h}(F)=i\,{\cal C}_h(F)-i\,\partial\bar\partial \log\eta
=i\,{\cal C}_h(F)-\eta^{-1}i\,\partial\bar\partial\eta+\eta^{-2}i\,\partial\eta\wedge \bar\partial\eta\,.
$$
The $L^2$-Bochner-Kodaira-Nakano Inequality whith respect to the hermitian metric $\eta h$ implies
\begin{eqnarray*}
&&\left([i\,{\cal C}_{\eta h}(F), \omega^*]v,v\right)_{\omega, \eta h}
\leq
\|\bar\partial_F v\|^2_{\omega,\eta h}+\|\bar\partial_{F,\eta}^*v\|^2_{\omega,\eta h}
\\
\\
&\leq& \|\bar\partial_F v\|^2_{\omega,\eta h}+\|\bar\partial_{F}^*v\|^2_{\omega,\eta h}+\|\eta^{-2}(\bar\partial \eta)^*v\|^2_{\omega,\eta h}
-2\Re e\left( \bar\partial_{F}^*v \,,\, (\bar\partial \eta)^*v\right)_{\omega, h}\,.
\end{eqnarray*}
By the identity \eqref{nor-cont} we infer
\begin{eqnarray*}
&&\left([\eta\,i\,{\cal C}_{ h}(F)-i\,\partial\bar\partial \eta\,, \omega^*]v,v\right)_{\omega, h}
\\
\\
&\leq& 
\|\bar\partial_F v\|^2_{\omega,\eta h}+\|\bar\partial_{F}^*v\|^2_{\omega,\eta h}
-2\Re e\left( \bar\partial_{F}^*v \,,\, (\bar\partial \eta)^*v\right)_{\omega, h}
\\
\\
&\leq&
\|\bar\partial_F v\|^2_{\omega,\eta h}+\|\bar\partial_{F}^*v\|^2_{\omega,(\eta +\lambda)h}
+
\|\lambda^{-2}(\bar\partial \eta)^*v\|^2_{\omega,\lambda h}
\,.
\end{eqnarray*}
Then the conclusion follows by applying again the identity \eqref{nor-cont}.\hfill$\Box$
\\
\\
Consider now the Hilbert spaces $H_{q+1}:=L^2_{\omega, h}(X,\Lambda^{n,q}T^*_X\otimes F)$, $q=0,1,2$ and the closed and densely defined linear operators 
$$
Tu:=\bar\partial_F\left((\eta+\lambda)^{1/2}\,u\right)\qquad\mbox{and}\qquad
Sv:=\eta^{1/2}\,\bar\partial_F v\,.
$$
We remark the identities $S\circ T=0$,  $T^*=(\eta+\lambda)^{1/2}\,\bar\partial^*_F$, $D(T^*)=D(\bar\partial^*_F)$ and $D(S)=D(\bar\partial_F)$. The Perturbed $L^2$-Bochner-Kodaira-Nakano Inequality \eqref{Prt-Geom-BKN} implies that the operators $T$ and $S$ satisfy the inequality \eqref{Abs-BKN} in the abstract existence theorem \eqref{Abs-Ex-Hil}. We infer by the abstract existence result \ref{Abs-Ex-Hil} B the following corollary.
\begin{corol}\label{OT-Pert-Ex}
In the setting of lemma  \ref{Prt-BKN}, if $g\in L^2_{\omega, h}(X,\Lambda^{n,1}T^*_X\otimes F)$, $\bar\partial_Fg=0$ satisfies the  inequality 
\begin{eqnarray}\label{OT-Pert-Ex-Est}
|(g,v)_{\omega, h}|^2\leq C\left((L_{\omega, h}^{\eta,\lambda}+\rho^2 Q)v,v\right)_{\omega, h}\,,
\end{eqnarray}
for all $v\in L^2_{\omega, h}(X,\Lambda^{n,1}T^*_X\otimes F)$, $($with $\rho\ge0$ a constant$)$ then there exist a solution $(u,h)$ of the perturbed $\bar\partial$-equation
$\bar\partial_F u+\rho PQh=g$, which satisfies the $L^2$-estimate
$$
n!\int_X(\eta+\lambda)^{-1}\,i^{n^2}u\wedge _h\bar u\,+ \, \int_X|h|^2_{\omega, h}\, dV_{\omega}\,\le \,C\,.
$$
\end{corol}
This is the case if $g$ is of the type $g=\bar\partial \eta\wedge \beta$, with $\beta\in L^2_h(X,K_X\otimes F)$ and 
\begin{eqnarray*}
\Theta^{\eta,\lambda}_h &\ge &\gamma\,i\partial \eta\wedge \bar\partial \eta\qquad\mbox{over an open set}\quad W\subset X\,,
\\
\\
\Theta^{\eta,\lambda}_h &\ge &\gamma\,i\partial \eta\wedge \bar\partial \eta-\rho^2\omega\qquad\mbox{over}\quad X\,,
\end{eqnarray*}
with $\gamma\ge 0$ over $X$ and $\gamma\ge k^{-1}$ over the support of $\beta$ for some constant $k>0$. In fact by lemma \ref{O-T-Key-Es} we deduce
\begin{eqnarray*}
|\left<g,v\right>_{\omega,h}|^2 &\le&  |\beta|^2_{\omega,h}\left<[i\partial \eta\wedge \bar\partial \eta, \omega^*]v,v\right>_{\omega,h}
\\
\\
&\le& k |\beta|^2_{\omega,h}\left<[\gamma i\partial \eta\wedge \bar\partial \eta, \omega^*]v,v\right>_{\omega,h}\,,
\end{eqnarray*}
over $X$.
We infer
\begin{eqnarray*}
\left|\int_W\left<g,v\right>_{\omega,h}dV_{\omega}\right|^2 &\le& C/2 \int_W\left<L_{\omega, h}^{\eta,\lambda}v,v\right>_{\omega,h}dV_{\omega}\,,
\\
\left|\int_{X\smallsetminus W}\left<g,v\right>_{\omega,h}dV_{\omega}\right|^2 &\le&  C/2\int_{X\smallsetminus W}\left<\left(L_{\omega, h}^{\eta,\lambda}+\rho^2I\right)v,v\right>_{\omega,h}dV_{\omega}\,,
\end{eqnarray*}
with $C=2k\,n!\int_Xi^{n^2}\beta\wedge _h\bar\beta$.
Thus the fundamental inequality \eqref{OT-Pert-Ex-Est} hold with $Q$ being the characteristic function of the set $X\smallsetminus W$.
\\
\\
{\bf Hermitian norms of forms.}
Let $(X,\omega )$ be a hermitian manifold. Let $h^*$ the corresponding hermitian metric over the complex vector bundle $T^*_{_{X,J}}$. In local complex coordinates we have the expressions
$$
h^*=4\,\sum_{k,l}\,\omega^{l\bar{k}}\,\frac{\partial }{\partial z_k}\otimes \frac{\partial }{\partial \bar{z}_l}.
$$
We remind that if $(V,J)$ is a complex vector space equipped with a hermitian metric $h$ then the corresponding hermitian metric $h_{_\C}$ over the complexified vector space $(V\otimes_{_\R}\C, i)$ is defined by the formula
$$
2h_{_\C}(v,w):=h(v,\overline{w})+\overline{h(\overline{v},w)},\quad v,w\in V\otimes_{_\R}\C,
$$
where we still note by  $h$ the $\C$-linear extension of $h$. 
Thus $h_{_\C}$ coincides with the sesquilinear extension over $V\otimes_{_\R}\C$ of the Riemann metric associated to $h$.
We infer the induced hermitian product on the vector bundle $\Lambda^{p,q}_{J}T^*_X$ is  given by the formula
\begin{eqnarray*}
&&\left<\Lambda_{j=1}^p\alpha_{1,j}\wedge \Lambda_{j=1}^q\beta_{1,j}\,,\,\Lambda_{j=1}^p\alpha_{2,j}\wedge \Lambda_{j=1}^q\beta_{2,j}\right>
\\
\\
&:=&(p+q)!\det\left(2^{-1}h^*(\alpha_{1,j},\bar\alpha_{2,l})\right)\,\overline{\det\left(2^{-1}h^*(\bar\beta_{1,j},\beta_{2,l})\right)}.
\end{eqnarray*}
In particular if $\alpha\in \Lambda^{p,0}_{J}T^*_X\otimes_{_{\C}}F$, with $(F,H)$ a hermitian vector bundle over $X$, then
$$
|\alpha|^2_{\omega, H}=\frac{n!\,p!}{(n-p)!}\,\frac{i^{p^2}\alpha\wedge _{_H}\bar \alpha\wedge \omega^{n-p}}{\omega^n}\,.
$$
\\
\\
We will need also the following lemma. (See \cite{Dem1} for more general statements.) 
\begin{lem}\label{Zot}
Let $(F,h)$ be a hermitian line bundle over a hermitian manifold $(X,\omega)$ of complex dimension $n$ and let 
$g\in \Gamma(X,\Lambda^{n,q}T^*_X\otimes L)$, $q\geq 1$ be a measurable section such that 
\begin{eqnarray*}
\int_X\gamma^{-1}|g|^2_{\omega, h}\,dV_{\omega}<+\infty\,,
\end{eqnarray*}
for some $\gamma\in C^0(X,[0,c])$. Let also $\Theta$ be a smooth $(1,1)$-form such that $\Theta\geq  \gamma\,\omega$ and set $L_{\hat\omega}:=[\Theta, \hat\omega^*]$ for any hermitian form $\hat\omega\ge \omega$. Then hold the inequalities
\begin{eqnarray}
|g|^2_{\hat\omega, h}\,dV_{\hat\omega}&\leq& |g|^2_{\omega, h}\,dV_{\omega}\,,\label{nrm-Comp-inq}
\\\nonumber
\\
|\left(g,v\right)_{\hat\omega,h}| ^2&\leq &\left(L_{\hat\omega}v,v\right)_{\hat\omega, h}\cdot \int_X(q \gamma)^{-1}|g|^2_{\omega, h}dV_{\omega}\,,\label{pnct-Finq}
\end{eqnarray}
for all $v\in L^2_{\hat\omega,h}(X,\Lambda^{n,q}T^*_X\otimes L)$.
\end{lem}
$Proof$.  Let $(\zeta_k)_k\subset C^{\infty}(U,\Lambda^{1,0}T_X)$ be a local frame such that
\begin{eqnarray*}
\hat\omega=i\sum_k  \zeta^*_k\wedge \bar\zeta^*_k\,,\qquad \omega=i\sum_k \lambda_k\, \zeta^*_k\wedge \bar\zeta^*_k\,,\qquad g=\sum_{|J|=q}  g_J\,\zeta^*\wedge \bar\zeta^*_J\,,
\end{eqnarray*}
with $0<\lambda_1\le...\le \lambda_n\le 1$. We set $\lambda_J:=\sum_{j\in J}\lambda_j$ and $\Lambda_J:=\prod_{j\in J}\lambda_j$. For all $j\in J$ hold $\lambda_j\ge \Lambda_J$, thus
\begin{eqnarray}\label{InLam}
\lambda_J\ge q \,\Lambda_J\,.
\end{eqnarray}
We remark now that $dV_{\hat\omega}=(\lambda_1\cdots \lambda_n)^{-1}dV_{\omega}$ and
$$
|g|^2_{\omega, h}=(n+q)!\sum_{|J|=q}(\lambda_1\cdots \lambda_n)^{-1}\Lambda_J^{-1}|g_J|^2\,.
$$
Thus
$$
|g|^2_{\omega, h}\,dV_{\omega}=(n+q)!\sum_{|J|=q}\Lambda_J^{-1}|g_J|^2\,dV_{\hat\omega}\ge |g|^2_{\hat\omega, h}\,dV_{\hat\omega}\,,
$$
which proves \eqref{nrm-Comp-inq}. The inequality \eqref{pnct-Finq} follows by combining the Cauchy-Schwarz inequality with the inequality 
\begin{eqnarray}
\left<L_{\hat\omega, h}^{-1}\,g,g\right>_{\hat\omega,h} dV_{\hat\omega}&\leq& (q \gamma)^{-1}|g|^2_{\omega, h}\,dV_{\omega}\,.\label{pnct-Comp-inq}
\end{eqnarray}
We prove now  \eqref{pnct-Comp-inq}. The conjugate of the identity \eqref{contrac-id} implies
\begin{eqnarray*}
[\omega, \hat\omega^*]g
\,=\,
\omega\wedge (\hat\omega^*\contract g)
&=&
\sum_{r,j=1}^n\sum_{|J|=q} \lambda_r\, g_J\,\zeta_r^*\wedge (\zeta_j\contract\zeta^*)\wedge \bar\zeta_r^*\wedge (\bar\zeta_j\contract\bar\zeta^*_J)
\\
\\
&=&
\sum_{j=1}^n\sum_{|J|=q} \lambda_j\, g_J\,\zeta^*\wedge\bar\zeta_j^*\wedge (\bar\zeta_j\contract\bar\zeta^*_J)\,.
\end{eqnarray*}
The fact that $\bar\zeta_j^*\wedge (\bar\zeta_j\contract\bar\zeta^*_J)=\bar\zeta^*_J$ if $j\in J$ and 
$\bar\zeta_j^*\wedge (\bar\zeta_j\contract\bar\zeta^*_J)=0$ if $j\not\in J$ implies
$$
[\omega, \hat\omega^*]g=\sum_{|J|=q} \lambda_J\, g_J\,\zeta^*\wedge \bar\zeta^*_J\,.
$$
The assumption $\Theta\ge \gamma\omega$ implies $L_{\hat\omega}\ge_{\,\hat\omega, h} \gamma[ \omega, \hat\omega^*]$ as operators, thus 
\begin{eqnarray}\label{InvCurv}
L_{\hat\omega}^{-1}\le_{\,\hat\omega, h} \gamma^{-1}[ \omega, \hat\omega^*]^{-1}\,.
\end{eqnarray}
The inequality \eqref{pnct-Comp-inq} will follow by combining \eqref{InvCurv} with the inequality
\begin{eqnarray*}
\left<[ \omega, \hat\omega^*]^{-1}g,g\right>_{\hat\omega,h} dV_{\hat\omega}&\leq& q ^{-1}|g|^2_{\omega, h}\,dV_{\omega}\,,
\end{eqnarray*}
that we prove now. In fact \eqref{InLam} implies
\begin{eqnarray*}
\left<[ \omega, \hat\omega^*]^{-1}g,g\right>_{\hat\omega,h} dV_{\hat\omega}
&=& 
(n+q)!\sum_{|J|=q}\lambda_J^{-1}|g_J|^2\,dV_{\hat\omega}
\\
\\
&=& (n+q)!\sum_{|J|=q}(\lambda_1\cdots \lambda_n)^{-1}\lambda_J^{-1}|g_J|^2\,dV_{\omega}
\\
\\
&\le&
q ^{-1}|g|^2_{\omega, h}\,dV_{\omega}\,.
\end{eqnarray*}
\hfill$\Box$
\begin{defini} A function $\psi\in C^0(X,\R_{\ge 0})$ over a topological space $X$ is called exhaustive if the open sets $X_c:=\{\psi<c\}$ are relatively compact for all $c>0$.
A complex manifold $X$ is called weakly pseudoconvex if there exist a smooth exhaustive function $\psi$ such that $i\partial\bar\partial \psi\ge 0$. 
\end{defini}
For domains $\Omega\subset \C^n$ the above weak pseudoconvexity notion is equivalent to pseudoconvexity. Note that every compact complex manifold is weakly pseudoconvex (take $\psi\equiv 0$). We observe also that the open sets $X_c$ are also weakly pseudoconvex. In fact the function
$$
\psi_c:=\log c-\log(c-\psi)\,,
$$ 
is exhaustive over $X_c$ and
$$
i\partial\bar\partial \psi_c=\frac{i\partial\bar\partial \psi}{c-\psi}+\frac{i\partial\psi\wedge \bar\partial \psi}{(c-\psi)^2}\ge 0\,,
$$
over $X_c$. We remind the following basic result \cite{Dem1}.
\begin{theorem}
Let $(X,\omega,\psi)$ be a weakly pseudoconvex K\"{a}hler manifold. Then the  K\"{a}hler metric $\omega+i\partial\bar\partial \psi^2$ is complete.
\end{theorem}
We will note by $B^p_{\delta}(0)\subset \C^p$ the $p$-times cartesian product of the disc $B_{\delta}(0)\subset \C$ with radius $\delta$ and center the origin.
\begin{defini}{\bf(Singularity condition)}
Let $Y\subset X$ be a pure $p$-dimensional complex analytic subset of a complex manifold  $X$. Let $\varphi\in L^1(X,\R)$ be a
quasi-plurisubharmonic function such that restriction $\varphi_{{|Y}}$ is  not identically $-\infty$ on any connected component of $Y$. 
We say that $\varphi$  is $Y$-admissible if the following condition hold.
\\
Let $U:=B^r_{\varepsilon}(0)\times B^p_{\delta}(0)\subset X$ be a coordinate open set with coordinates $(z,\zeta)$ centered in an arbitrary regular point $y\in Y$ such that $Y\cap U=\{z=0\}$ is smooth. Set 
$$
Y_z:=\{(z,\zeta)\,\mid\,\zeta\in B^p_{\delta}(0)\}\,,
$$ 
and let ${\cal J}({\varphi})\subset {\cal O}_X$ be the multiplier ideal sheaf associated to $\varphi$.
\\
For all $f\in H^0(U,{\cal J}({\varphi}))$ such that 
\begin{eqnarray}\label{mtp-Id-shf-L^2}
\int_{\zeta\in Y_0}|f|^2e^{-\varphi}\,i^{p^2}d\zeta\wedge d\bar\zeta<+\infty\,,
\end{eqnarray}
there exist a zero measure set $E\subset B^r_{\varepsilon}(0)\smallsetminus \{0\}$ and $\varepsilon'\in (0,\varepsilon)$ such that the map
$$
z\in B^r_{\varepsilon'}(0)\smallsetminus E\longmapsto \int_{\zeta\in Y_z}|f|^2e^{-\varphi}\,i^{p^2}d\zeta\wedge d\bar\zeta<+\infty\,,
$$
is continuous at the origin.
\end{defini}
The $Y$-admissibility condition is always satisfied for quasi-plurisubharmonic functions $\varphi$ which can locally be expressed as 
$$
\varphi=c\log\sum_j|h_j|^2+u\,,
$$
with $c\in \R_{>0}$, $h_j$ holomorphic, $u$ continuous (this functions are called with complex analytic singularities)
and such that restriction $\varphi_{{|Y}}$ is  not identically $-\infty$ on any connected component of $Y$. In fact in the case 
$$
\varphi^{-1}(-\infty)_{\varepsilon,\delta}:=\varphi^{-1}(-\infty)\cap \left(B^r_{\varepsilon}(0)\times B^p_{\delta}(0)\right)\subset Y_0\,,
$$
the $Y$-admissibility follows from the dominated convergence theorem, which can be applied thanks to the $L^2$-assumption \eqref{mtp-Id-shf-L^2}.
\\
In the case $\varphi^{-1}(-\infty)_{\varepsilon,\delta}\not\subset Y_0$ 
let $f\in H^0(U,{\cal J}({\varphi}))$ arbitrary. The assumption 
$$
Y_0\not\subset \varphi^{-1}(-\infty)_{\varepsilon,\delta}\,,
$$ 
implies the existence of a blow-up map $\mu:(z,\xi)\mapsto(z,\zeta)$ such that $f\circ \mu\simeq z^{\alpha}\xi^{\beta}$, up to an invertible factor and
$$
\varphi\circ \mu=c\log|\xi^{\gamma}|^2+R\,, 
$$
with $R$ continuous. This last equality follows from the fact that we can construct the blow-up map  in a way that 
$$
\mu^*\sum_j{\cal O}\cdot h_j
$$ 
is an invertible sheaf. 
Moreover we can also assume that the Jacobian $J(\mu)$ of $\mu$ equal to
a monomial $\xi^\delta$, up to an invertible factor. We infer
$$
z\longmapsto\int_{\zeta\in Y_z}|f|^2e^{-\varphi}\,i^{p^2}d\zeta\wedge d\bar\zeta \;=\;\int_{\xi\in B^p_{\delta}(0)}g(z,\xi)\,\frac{|z^{\alpha}|^2\;|\xi^{\beta}|^2}{\;|\xi^{\gamma}|^{2c}\;}\,|\xi^{\delta}|^2\,i^{p^2}d\xi\wedge d\bar\xi\,,
$$
with $g$ continuous and invertible.
\section{The geometric singular $L^2$-extension result.}\label{sng-ext-rst}
\begin{theorem}\label{Ext-tm}
Let $(L,h)$ and $(E,H)$ be two holomorphic hermitian vector bundles of rank $\rk_{_{\C}}L=1$, $\rk_{_{\C}}E=r$ over a complex manifold $X$ of complex dimension $n$ admitting a complex analytic subset $A$ such that $X\smallsetminus A$ is Stein, let $\sigma\in H^0(X,E)$ such that $|\sigma|_{_H}\le e^{-\alpha}$, $\alpha\in \R_{\ge 1}$ on $X$ and 
$$
Y:=\overline{Y}_0\not \subset A\,,\;\;Y_0:=\left\{x\in X\,|\,\sigma(x)=0\,,\Lambda^rd\sigma(x)\not=0\right\}\,.
$$
Consider also a quasi-plurisubharmonic function
$\varphi\in L^1(X,\R)$ such that  
the restriction $\varphi_{{|Y}}$ is  not identically $-\infty$ on any connected component of $Y$. Assume that $\varphi$  is $Y$-admissible and the curvature current
\begin{eqnarray*}
\Theta:=i\,{\cal C}_h(L)+i\partial\bar\partial \left(\varphi+r\log|\sigma|^2_{_H}\right)
\end{eqnarray*}
satisfies the positivity assumptions
\begin{eqnarray}\label{Cruc-Pos-hyp}
\Theta\ge 0
\quad\mbox{and}\quad
\Theta\ge \alpha^{-1}H\left(i{\cal C}_H(E)\sigma,\sigma\right)|\sigma|^{-2}_{_H}\,.
\end{eqnarray}
Then there exist a uniform constant $C_r>0$ depending only on $r$ and on a fixed cut off function
such that for any section $f\in H^0(Y,K_X\otimes L)$, with the $L^2$-property
\begin{eqnarray}
I_Y(f,\sigma,\varphi):=\int_Yi^{p^2}f_{\sigma}\wedge_{_{h,H}}\bar f_{\sigma}\;e^{-\varphi}<+\infty\,,
\end{eqnarray}
where $p:=\dim_{_{\C}}Y=n-r$ and 
$$
f_{\sigma}:=f/(\Lambda^rd\sigma)\in H^0(Y,K_Y\otimes L_{{|Y}}\otimes \Lambda^rE^*_{{|Y}})\,,
$$
there exist
$F\in H^0(X,K_X\otimes L)$ such that $F=f$ over $Y$ and
\begin{eqnarray*}
J_X(F,\sigma,\varphi):=\int_X\,\frac{\;i^{n^2}F\wedge_{_{h}} \bar F\;e^{-\varphi}\; }{\;|\sigma|^{2r}_{_H}\,(\log|\sigma|_{_H})^2\;}\,\le \,C_rI_Y(f,\sigma,\varphi)\,.
\end{eqnarray*}
\end{theorem}
We remark in particular the trivial inequalities
$$
e^{2\alpha(r-\varepsilon)}\int_Xi^{n^2}F\wedge_{_{h}} \bar F\;e^{-\varphi}\le \int_X\,\frac{\;i^{n^2}F\wedge_{_{h}} \bar F\;e^{-\varphi}\; }{\;|\sigma|^{2(r-\varepsilon)}_{_H}\;}
\le \varepsilon^{-2}J_X(F,\sigma,\varphi)
$$
for all $\varepsilon\in (0,r)$. We observe also that the condition $Y\not \subset A$ is always satisfied in the case $X$ is a complex projective variety. Moreover in the case $X$ projective (or bounded pseudoconvex domain) and $\rk_{_{\C}}E=1$ the $Y$-admissibility condition of $\varphi$ is not needed (see subsection \ref{Loc-Proj}). 

In particular  geometric situations the curvature conditions \eqref{Cruc-Pos-hyp} can be drastically simplified and the $Y$-admissibility condition can be dropped (see the remark at the end of the subsection \ref{Loc-Proj}). 
\\
\\
{\bf Proof.} The proof is divided in several steps.
\\
\\
{\bf A). The positivity estimate.} We set 
$$
S_{\varepsilon}:=\log(|\sigma|^2+\varepsilon^2)\,.
$$
Then a quite standard computation \cite{De-Pa} implies
\begin{eqnarray}\label{Hol-Pos-est}
i\partial\bar\partial S _{\varepsilon}
&\ge&
\frac{\varepsilon^2}{|\sigma|^2}\,i\partial S _{\varepsilon}\wedge \bar\partial S _{\varepsilon}-\Phi_{\varepsilon}\,,
\quad\mbox{with}\quad
\Phi_{\varepsilon}:=\frac{H(i{\cal C}_H(E)\sigma,\sigma)}{|\sigma|^2+\varepsilon^2}\,.\qquad
\end{eqnarray}
On the other hand hold the identity
\begin{eqnarray}\label{Lib-Pos-EstI}
i\partial\bar\partial S _{\varepsilon}=\frac{|\sigma|^2}{|\sigma|^2+\varepsilon^2} \,i\partial\bar\partial\log|\sigma|^2+\frac{\varepsilon^2}{|\sigma|^2}\,i\partial S_{\varepsilon}\wedge \bar\partial S _{\varepsilon}\,.
\end{eqnarray}
We consider also the family of locally bounded functions 
$$
\varphi_{\delta}:=\log (e^{\varphi}+\delta^2)
$$
which decreases to $\varphi$ as $\delta\rightarrow 0$ and with complex hessian
\begin{eqnarray}\label{Lib-Pos-EstII}
i\partial\bar\partial \varphi_{\delta}=\frac{e^{\varphi}}{e^{\varphi}+\delta^2}\,i\partial\bar\partial \varphi +\frac{\delta^2}{e^{\varphi}} \,i\partial\varphi_{\delta}\wedge \bar\partial \varphi_{\delta}\,.
\end{eqnarray}
We introduce the following notation. Let $U, V\subset X$ be two sets. We will note  $U_V:=U\smallsetminus V$.
Let $\Sigma:=Y_{Y_0}$. We observe as in \cite{Dem3} that the fact that $X_A:=X\smallsetminus A$ is Stein implies the existence of a complex hypersurface $Z\subset X_A$ such that $\Sigma_A\subset Z$ and $Y_A\not\subset Z$. 

Moreover if we set $A':=A\cup Z$ then 
$X_{A'}=X_A\smallsetminus Z$ is also Stein and $Y_{A'}$ is a smooth non empty subvariety of $X_{A'}$. 
So let $\psi\in C^{\infty}(X_{A'},\R_{\ge 0})$ be an exhaustive function such that $\omega:=i\partial\bar\partial \psi>0$. We will consider also the K\"{a}hler manifold $X_c:=\{\psi<c\}\subset\subset X_{A'}$ equipped with the complete K\"{a}hler metric $\omega_c:=\omega+i\partial\bar\partial \psi^2_c$.
For all $\delta, \tau>0$ we introduce the current over $X_{A'}$
$$
\Theta_{\delta,\tau}:=i\,{\cal C}_h(L)+i\partial\bar\partial(\varphi_{\delta}+\tau\psi+rS_{\delta})\,.
$$
Let $K>1$ be a sufficiently big constant such that 
$$
i\,{\cal C}_h(L)+i\partial\bar\partial \varphi_{\delta}\ge -K\omega\,,
$$  
$$
i\,{\cal C}_h(L)+r i\partial\bar\partial\log|\sigma|^2\geq -K\omega\,,
$$ 
and $\Theta_{\delta,\tau}\ge -K\,\omega_c$ over $X_c$. 
We infer by \eqref{Lib-Pos-EstI}, \eqref{Lib-Pos-EstII} and \eqref{Hol-Pos-est} the inequalities over $X_{A'}$
\begin{eqnarray*}
\Theta_{\delta,\tau}&\geq&i\,{\cal C}_h(L)+i\partial\bar\partial (\varphi_{\delta}+\tau\psi)+\frac{|\sigma|^2}{|\sigma|^2+\delta^2} \,ri\partial\bar\partial\log|\sigma|^2
\\
\\
&=&
\tau\omega+\frac{\delta^2}{|\sigma|^2+\delta^2}\left(i\,{\cal C}_h(L)+i\partial\bar\partial \varphi_{\delta}\right)
\\
\\
&+&
\frac{|\sigma|^2}{|\sigma|^2+\delta^2}\left[i\,{\cal C}_h(L)+i\partial\bar\partial\left(\varphi_{\delta}+r\log|\sigma|^2\right)\right]
\\
\\
&\ge&
\left(\tau-\frac{\delta^2K}{|\sigma|^2+\delta^2}\right)\omega
\\
\\
&+&
\frac{|\sigma|^2}{|\sigma|^2+\delta^2}\left(i\,{\cal C}_h(L)+ri\partial\bar\partial\log|\sigma|^2+\frac{e^{\varphi}}{e^{\varphi}+\delta^2}\,i\partial\bar\partial \varphi\right)
\\
\\
&\ge&
\left(\tau-\frac{\delta^2K}{|\sigma|^2+\delta^2}\right)\omega
+
\frac{|\sigma|^2}{|\sigma|^2+\delta^2}\,\frac{\delta^2}{e^{\varphi}+\delta^2}\,\left(i\,{\cal C}_h(L)+ri\partial\bar\partial\log|\sigma|^2\right)
\\
\\
&+&
\frac{|\sigma|^2}{|\sigma|^2+\delta^2}\,\frac{e^{\varphi}}{e^{\varphi}+\delta^2}\,\Theta
\\
\\
&\ge&
\tau_{\delta}\,\omega+\frac{|\sigma|^2}{|\sigma|^2+\delta^2}\,\frac{e^{\varphi}}{e^{\varphi}+\delta^2}\,\Theta\,,
\end{eqnarray*}
with 
$$
\tau_{\delta}:=\tau-\frac{\delta^2K}{|\sigma|^2+\delta^2}-\frac{|\sigma|^2}{|\sigma|^2+\delta^2}\,\frac{\delta^2K}{e^{\varphi}+\delta^2}\,.
$$
The hypothesis $\Theta\ge 0$ implies 
\begin{eqnarray}\label{Key-PosI}
\Theta_{\delta,\tau}\ge \tau_{\delta}\,\omega+\frac{|\sigma|^2}{|\sigma|^2+\varepsilon^2}\,\frac{e^{\varphi}}{e^{\varphi}+\delta^2}\,\Theta\,.
\end{eqnarray}
for all $\delta\in (0,\varepsilon)$.
We introduce now the function $\chi:(-\infty, 0]\rightarrow (-\infty, 0]$, 
$$
\chi(t):=t-\log(1-t)\,,
$$
and we observe the trivial inequalities $\chi\le t$, $1\le \chi'\le 2$, $\chi''=(1-t)^{-2}$.
We set $\eta_{\varepsilon}:=\varepsilon-\chi(S _{\varepsilon})$ and we remark the inequality $\varepsilon-S_{\varepsilon}\le \eta_{\varepsilon}\le \varepsilon-2S_{\varepsilon}$. 
In particular 
$$
\eta_{\varepsilon}\ge 2\alpha \ge \alpha\chi'\,,
$$ 
for $\varepsilon>0$ sufficiently small. 
By deriving the identity  $\partial\eta_{\varepsilon}=-\chi'(S_{\varepsilon})\partial S_{\varepsilon}$ we infer the inequalities over $X$
\begin{eqnarray*}
-i\partial\bar\partial\eta_{\varepsilon}
&=&
\chi''(S _{\varepsilon})\,i\partial S _{\varepsilon}\wedge \bar\partial S _{\varepsilon}\,+\,\chi'(S _{\varepsilon})\,i\partial\bar\partial S _{\varepsilon}
\\
\\
&\ge &
\left(\, \chi''(S _{\varepsilon})\,+\,\, \frac{\varepsilon^2}{|\sigma|^2}\,\chi'(S _{\varepsilon})  \,\right)\,i\partial S _{\varepsilon}\wedge \bar\partial S _{\varepsilon}
\,-\,\chi'(S _{\varepsilon})\Phi_{\varepsilon}
\\
\\
&\ge &
\left(\, \frac{\chi''(S _{\varepsilon})}{\chi'(S _{\varepsilon})^2} \,+\, \frac{\varepsilon^2}{|\sigma|^2\,\chi'(S _{\varepsilon})}  \,\right)\,i\partial \eta _{\varepsilon}\wedge \bar\partial \eta _{\varepsilon}
\,-\,\chi'(S _{\varepsilon})\Phi_{\varepsilon}\,.
\end{eqnarray*}
Thus if we set $\lambda_{\varepsilon}:=\chi'(S _{\varepsilon})^2/\chi''(S _{\varepsilon})$ we infer the inequality
\begin{eqnarray}\label{Part-Pos-Est}
-i\partial\bar\partial\eta_{\varepsilon}-\lambda^{-1}_{\varepsilon}\,i\partial \eta _{\varepsilon}\wedge \bar\partial \eta _{\varepsilon}
\,\ge \,
\frac{\varepsilon^2}{2|\sigma|^2}\,i\partial \eta _{\varepsilon}\wedge \bar\partial \eta _{\varepsilon}
\,-\,\chi'(S _{\varepsilon})\Phi_{\varepsilon}\,.
\end{eqnarray}
Combining  \eqref{Key-PosI} with the positivity assumptions \eqref{Cruc-Pos-hyp} and with the inequality $\eta_{\varepsilon}\ge \alpha\chi'$, we infer 
\begin{eqnarray*}
\eta_{\varepsilon}\Theta_{\delta,\tau}
\;\ge\;
\eta_{\varepsilon}\tau_{\delta}\,\omega\,+\,
\frac{|\sigma|^2}{|\sigma|^2+\varepsilon^2}\,\frac{e^{\varphi}\,\alpha\,\chi'(S _{\varepsilon})}{e^{\varphi}+\delta^2}\,\Theta
\;\ge\;
\eta_{\varepsilon}\tau_{\delta}\,\omega\,+\,
\frac{e^{\varphi}\chi'(S _{\varepsilon})}{e^{\varphi}+\delta^2}\,\Phi_{\varepsilon}
\,.
\end{eqnarray*}
for all $\delta\in (0,\varepsilon)$.
This combined with the inequality \eqref{Part-Pos-Est} yields
\begin{eqnarray}\label{Cruc-Pos-Est}
\Theta^{\varepsilon}_{\delta,\tau}
&:=& 
\eta_{\varepsilon}\Theta_{\delta,\tau}-i\partial\bar\partial\eta_{\varepsilon}-\lambda_{\varepsilon}^{-1}i\,\partial\eta_{\varepsilon}\wedge \bar\partial\eta_{\varepsilon}\nonumber
\\\nonumber
\\
&\ge&
\frac{\varepsilon^2}{2|\sigma|^2}\,i\,\partial\eta_{\varepsilon}\wedge \bar\partial\eta_{\varepsilon}
+\eta_{\varepsilon}\tau_{\delta}\,\omega-\frac{\delta^2\chi'(S _{\varepsilon})}{e^{\varphi}+\delta^2}\,\Phi_{\varepsilon}\nonumber
\\\nonumber
\\
&\ge&
\frac{\varepsilon^2}{2|\sigma|^2}\,i\,\partial\eta_{\varepsilon}\wedge \bar\partial\eta_{\varepsilon}
+\tau_{\varepsilon,\delta}\,\omega\,,\quad
\end{eqnarray}
with
$$
\tau_{\varepsilon,\delta}:=2\tau-\frac{\delta^2K_{\varepsilon}}{|\sigma|^2+\delta^2}-\frac{\delta^2K_{\varepsilon}}{e^{\varphi}+\delta^2}-\frac{|\sigma|^2}{|\sigma|^2+\delta^2}\,\frac{\delta^2K_{\varepsilon}}{e^{\varphi}+\delta^2}\,,
$$
and $K_{\varepsilon}:=K(\varepsilon-4\log\varepsilon)\ge K>1$, for $\varepsilon>0$ sufficiently small.
\\
\\
{\bf The lost of positivity locus.}
\\
Let $\tau\in (0,1)$ and set $V_{\delta}:=X_c\cap\{\tau_{\varepsilon,\delta}<\tau\}$. An elementary computation shows that $V_{\delta}=X_c\cap (E'_{\delta}\cup  E''_{\delta})$, with
\begin{eqnarray*}
E'_{\delta}&:=&\left\{x\in X\,\mid \,|\sigma(x)|^2\le \delta^2(K_{\varepsilon}/\tau-1)\right\}\,,
\\
\\
E''_{\delta}&:=&\left\{x\in X\smallsetminus E'_{\delta}\,\mid\,\varphi(x)<2\log \delta+\log (2K_{\varepsilon}-\tau)+\log M_{\delta}(x)\right\}\,,
\\
\\
M_{\delta}&:=&\frac{|\sigma|^2+\delta^2}{\tau|\sigma|^2+(\tau-K_{\varepsilon})\delta^2}\,.
\end{eqnarray*}
We infer the family $(V_{\delta})_{\delta>0}$ is non-decreasing and 
\begin{eqnarray}\label{IntesYY}
X_c\cap \left(Y\cup \varphi^{-1}(-\infty)\right)=\bigcap_{\delta>0}V_{\delta}\,.
\end{eqnarray}
Moreover the fact that the function $\varphi$ is upper semicontinuous implies that the set $E''_{\delta}$ is open. Let $U_{\delta}\subset X_c$ be the interior of $X'_c\smallsetminus V_{\delta}$.
By elementary facts about measure theory we infer the identity
\begin{eqnarray}\label{Zer-Meas}
\int_{(X_c\smallsetminus V_{\delta})\smallsetminus U_{\delta}}\omega^n_c=0\,,
\end{eqnarray}
which will be very useful in step (C). We denote by $\Theta^{\varepsilon}_{\delta,\tau,t}$ the current obtained by formally replacing $\varphi_{\delta}$ with $\varphi_{\delta,t}$ in the definition of $\Theta^{\varepsilon}_{\delta,\tau}$.
By \eqref{Cruc-Pos-Est} and \cite{Dem1}, we infer the existence of a regularising family of locally uniformly bounded from above (over $X$) smooth functions
$(\varphi_{\delta, t})_{t>0}\subset C^{\infty}(X,\R)$ such that $\varphi_{\delta, t}\downarrow \varphi_{\delta}$ as $t\rightarrow 0$ and 
\begin{eqnarray}\label{G-Lost-Pos}
\Theta^{\varepsilon}_{\delta,\tau,t}
\ge
\frac{\varepsilon^2}{2|\sigma|^2}\,i\,\partial\eta_{\varepsilon}\wedge \bar\partial\eta_{\varepsilon}
+(\tau-\mu_t)\,\omega\,,
\end{eqnarray}
over $U_{\delta}$,
with $(\mu_t)_{t>0}\equiv(\mu_t^{\varepsilon,\delta,\tau})_{t>0}\subset C^0(U_{\delta}, \R_{>0})$ such that $\mu_t\downarrow 0$ locally uniformly as $t\rightarrow 0$. 
On the other hand by the definition of $\Theta^{\varepsilon}_{\delta,\tau,t}$ and \eqref{Part-Pos-Est} we infer the inequality
\begin{eqnarray}\label{MRG-Lost-Pos}
\Theta^{\varepsilon}_{\delta,\tau,t}\ge \frac{\varepsilon^2}{2|\sigma|^2}\,i\,\partial\eta_{\varepsilon}\wedge \bar\partial\eta_{\varepsilon}-k^2_{c,\varepsilon}\,\omega_c\,,
\end{eqnarray}
over $X_c$, with $k_{c,\varepsilon}>0$ a constant uniform in the parameters $\delta, \tau$ and $t$.
\\
\\
{\bf Step B). Construction of a $L^2$-extension over $X_c$.} 
\\
Let $(B_j)_j$ be a finite family of coordinate open balls $B_j\subset X_c$ which covers $Y_c:=X_c\cap Y$ such that $B_j\cap Y$ is a vector space with respect to the complex coordinates of $B_j$ and the bundles $K_X$ and $L$ are holomorphically trivial over $B_j$. Let $\tilde f_j$ be a holomorphic extension of $f$ over $B_j$.
By using a partition of unity $(\theta_j)_j$ subordinated to $(B_j)_j$ one can construct a section 
$$
F_{\infty}:=\sum_j\theta_j\,\tilde f_j\in C^{\infty}(X_c,K_X\otimes L)\,,
$$ 
such that $F_{\infty}=f$ and $\bar{\partial}_LF_{\infty}=0$ over $Y$. 
In fact 
$$
\bar{\partial}_LF_{\infty}=\sum_j\bar{\partial}\theta_j\wedge \tilde f_j\,,
\quad\mbox{thus}
\quad
\bar{\partial}_LF_{\infty}=\sum_j\bar{\partial}\theta_j\wedge  f=0\quad  \mbox{over}\quad Y\,.
$$
Consider now  $\theta\in C^{\infty}(\R,[0,1])$ such that $\theta=1$ over $(-\infty, 1/2]$, $\theta=0$ over $[1,+\infty)$ and $|\theta'|\le 3$. We set $F_{\varepsilon}:=\theta(\varepsilon^{-2}|\sigma|^2)F_{\infty}$. 
We will solve the equation
\begin{eqnarray}\label{O-T-eq}
\left  \{
\begin{array}{lr}
\bar{\partial}_LU_{\varepsilon}\,=\, \bar{\partial}_LF_{\varepsilon}
\\
\\
U_{\varepsilon\,{|Y}}\,=\,0  .
\end{array}
\right.
\end{eqnarray}
In order to get $U_{\varepsilon\,{|Y}}\,=\,0$ we just need to inshur that the function $|U_{\varepsilon}|^2|\sigma|^{-2r}$ is locally integrable near $Y$ since the complex codimension of $Y$ is $r$. In this way we will obtain a holomorphic section $H_{\varepsilon}:=F_{\varepsilon}-U_{\varepsilon}$ over $X_c$, which coincides with $f$ over $Y_c$.
\\
The definition of $\Theta^{\varepsilon}_{\delta,\tau,t}$ combined with the inequality \eqref{Part-Pos-Est} implies the existence of a sufficiently big $\tau>1$ such that 
\begin{eqnarray}\label{Rf-Pos-Es}
\Theta^{\varepsilon}_{\delta,\tau,t}\ge \frac{\varepsilon^2}{2|\sigma|^2}\,i\,\partial\eta_{\varepsilon}\wedge \bar\partial\eta_{\varepsilon} +\omega\,,
\end{eqnarray}
over $X_c$ for all $\varepsilon,\delta,t\in (0,1)$. We set 
$
h'\equiv h_{\delta,\tau,t}:=(|\sigma|^2+\delta^2)^{-r}h\,e^{-\varphi_{\delta,t}-\tau\psi}
$.
By using the identity $\varepsilon^{-2}e^{S_{\varepsilon}}=1+\varepsilon^{-2}|\sigma|^2$ as in \cite{Dem3}, one can decompose $g_{\varepsilon}:=\bar\partial_LF_{\varepsilon}=g_{1,\varepsilon}+g_{2,\varepsilon}$, with
\begin{eqnarray*}
g_{1,\varepsilon} &:=&
(1+\varepsilon^{-2}|\sigma|^2)\,\theta'(\varepsilon^{-2}|\sigma|^2)\,\partial S_{\varepsilon}\wedge F_{\infty}
\\
\\
&=&
(1+\varepsilon^{-2}|\sigma|^2)\,\theta'(\varepsilon^{-2}|\sigma|^2)\,\chi'(S _{\varepsilon})^{-1}\,\partial \eta_{\varepsilon}\wedge F_{\infty}\,,
\\
\\
g_{2,\varepsilon}&:=& \theta(\varepsilon^{-2}|\sigma|^2)\,\partial_LF_{\infty}\,.
\end{eqnarray*}
We remark that the support of $g_{j,\varepsilon}$ is contained in $\{|\sigma|<\varepsilon\}$, thus over this set hold the trivial inequalities $2\varepsilon^{-2}|\sigma|^2\le 2$ and $(1+\varepsilon^{-2}|\sigma|^2)\,\chi'(S _{\varepsilon})^{-1}\,\le\,2$. So by combing the lemma \ref{O-T-Key-Es} with the inequality \ref{Rf-Pos-Es} we infer the following punctual estimate
\begin{eqnarray}\label{Fund-Ex-es}
\left|\left<g_{1,\varepsilon}\,,v\right>_{\omega_c,h'}\right|^2
&\le&
4\theta'(\varepsilon^{-2}|\sigma|^2)^2\,|F_{\infty}|^2_{\omega_c,h'}\left<[i\,\partial\eta_{\varepsilon}\wedge \bar\partial\eta_{\varepsilon}\,,\,\omega_c^*]v,v\right>_{\omega_c,h'}\nonumber
\\\nonumber
\\
&\le&
8\theta'(\varepsilon^{-2}|\sigma|^2)^2\,|F_{\infty}|^2_{\omega_c,h'}\left<\left[\frac{\varepsilon^2}{2|\sigma|^2}\,i\,\partial\eta_{\varepsilon}\wedge \bar\partial\eta_{\varepsilon}\,,\,\omega_c^*\right]v,v\right>_{\omega_c,h'}\nonumber
\\\nonumber
\\
&\le&
8\theta'(\varepsilon^{-2}|\sigma|^2)^2\,|F_{\infty}|^2_{\omega_c,h'}\left<L_{\varepsilon}v,v\right>_{\omega_c,h'}\,,
\end{eqnarray}
with $L_{\varepsilon}:=[\Theta^{\varepsilon}_{\rho,\delta,\tau,t},\omega_c^*]$. By lemma \ref{Zot} we infer the estimate
\begin{eqnarray*}
\left|\left(g_{2,\varepsilon}\,,v\right)_{\omega_c,h'}\right|^2\le \left(L_{\varepsilon}v,v\right)_{\omega_c,h'}\int_{X_c\cap\{|\sigma|<\varepsilon\}}|\bar\partial F_{\infty}|^2_{\omega,h}\,|\sigma|^{-2r}e^{-\varphi_{\delta,t}-\tau\psi}\,dV_{\omega}\,.
\end{eqnarray*}
The Taylor expansion of $\bar\partial_LF_{\infty}$ near $Y_c$ implies $|\bar\partial_LF_{\infty}|^2_{\omega,h}\le C|\sigma|^2$ over the set $X_c\cap\{|\sigma|<\varepsilon\}$.
By putting all this together we infer for all $v\in L^2_{\omega_c,h'}(X_c,\Lambda^{n,1}T^*_X\otimes L)$, the estimate
\begin{eqnarray*}
\left|\left(g_{\varepsilon},v\right)_{\omega_c,h'}\right|^2 \le 2C'_{\varepsilon}\left(L_{\varepsilon}v,v\right)_{\omega_c,h'}\,,
\end{eqnarray*}
with $C'_{\varepsilon}:=C_{\varepsilon}+c_{\varepsilon}$ and
\begin{eqnarray*}
C_{\varepsilon}&:=&8n!\int_{X_c}\theta'(\varepsilon^{-2}|\sigma|^2)^2\,i^{n^2}F_{\infty}\wedge_{h} \overline{F}_{\infty}\,|\sigma|^{-2r}e^{-\varphi_{\delta,t}-\tau\psi}\,,
\\
\\
c_{\varepsilon}&:=&C\int_{X_c\cap\{|\sigma|<\varepsilon\}}|\sigma|^{-2(r-1)}e^{-\varphi_{\delta,t}-\tau\psi}dV_{\omega}\longrightarrow 0\,,\quad\mbox{as}\quad \varepsilon\rightarrow 0\,.
\end{eqnarray*}
By corollary \ref{OT-Pert-Ex} we infer the existence of a solution of the $\bar\partial$-problem \eqref{O-T-eq} with the $L^2$-estimate
$$
n!\int_{X_c}(\eta_{\varepsilon}+\lambda_{\varepsilon})^{-1}\,i^{n^2}U_{\varepsilon} \wedge_h \overline{U}_{\varepsilon}\;|\sigma|^{-2r}e^{-\varphi_{\delta,t}-\tau\psi}\le 2C'_{\varepsilon}\,.
$$
We need now the following estimate that we will prove in the step (D).
\begin{eqnarray}\label{Ban-Est}
\eta_{\varepsilon}+\lambda_{\varepsilon}\le \left(5+O(\varepsilon)\right)S_{\varepsilon}^2
\end{eqnarray}
Combining this with the previous inequality we infer
\begin{eqnarray*}
n!\int_{X_c}S_{\varepsilon}^{-2}\,i^{n^2}U_{\varepsilon} \wedge_h \overline{U}_{\varepsilon}\;(|\sigma|^2+\varepsilon^2)^{-r}\,e^{-\varphi_{\delta,t}-\tau\psi}\le \left(5+O(\varepsilon)\right)2C'_{\varepsilon}\,.
\end{eqnarray*}
Moreover the fact that the section $F_{\varepsilon}$ has uniformly bounded norm and support contained in $X_c\cap\{|\sigma|<\varepsilon\}$ implies the estimate
\begin{eqnarray*}
&&n!\int_{X_c}S_{\varepsilon}^{-2}\,i^{n^2}F_{\varepsilon}\wedge \overline{F}_{\varepsilon}\,(|\sigma|^2+\varepsilon^2)^{-r}\,e^{-\varphi_{\delta,t}-\tau\psi}
\\
\\
&\le&
C(\log 2\varepsilon^2)^{-2}\int_{X_c\cap \{|\sigma|<\varepsilon\}} \varepsilon^{-2r}dV
\;\le\;
C'(\log 2\varepsilon^2)^{-2}\,,
\end{eqnarray*}
for some constants $C,C'>0$ uniform in the parameters $\varepsilon,\delta,t$.
Combining this with the previous inequality we infer the last of the following inequalities
\begin{eqnarray}
&&
\alpha_{\mu}n!\int_{X_c}\,i^{n^2}H_{\varepsilon}\wedge \overline{H}_{\varepsilon}\,e^{-\varphi_{\delta,t}-\tau\psi}\nonumber
\\\nonumber
\\
&\le&
\mu^2n!\int_{X_c}\,i^{n^2}H_{\varepsilon}\wedge \overline{H}_{\varepsilon}\,(|\sigma|^2+\varepsilon^2)^{-(r-2\mu)}\,e^{-\varphi_{\delta,t}-\tau\psi}\nonumber
\\\nonumber
\\
&\le&
n!\int_{X_c}S_{\varepsilon}^{-2}\,i^{n^2}H_{\varepsilon}\wedge \overline{H}_{\varepsilon}\,(|\sigma|^2+\varepsilon^2)^{-r}\,e^{-\varphi_{\delta,t}-\tau\psi}\nonumber
\\\nonumber
\\
&\le&
2\left(5+O(\varepsilon)\right) C'_{\varepsilon}+2C'(\log 2\varepsilon^2)^{-2}\,,\label{tc-L2-OT}
\end{eqnarray}
where $\alpha_{\mu}>0$ is a constant depending only on $\mu\in (0,r/2)$ and $\alpha$. By the lemma \ref{Conv-Int-lm}, that we will prove in sequel, we infer the convergence
\begin{eqnarray}\label{G-conv-Int-Y}
\lim_{\varepsilon\rightarrow 0} C_{\varepsilon}=n!8\,K_r\,I_{Y_c}(f,\sigma,\varphi_{\delta,t}+\tau\psi)<+\infty
\end{eqnarray}
By weak compactness we can extract a weakly convergent sequence $H_{\varepsilon_k}$ with limit $H_{\delta,t}$ as $\varepsilon_k\rightarrow 0$. 
Moreover  $H_{\delta,t}=f$ over $Y_c$ as follows directly from the definition of Bergman Kernel.
The fact that $(-S_{\varepsilon})^{-1}$ is uniformly bounded implies by claim \ref{Prod-w-lim}
$$
(-S_{\varepsilon_k})^{-1}H_{\varepsilon_k}\longrightarrow (-2\log|\sigma|)^{-1}H_{\delta,t}
$$ 
weakly as $\varepsilon_k\rightarrow 0$. By claim \ref{w-lim-wei} (A), \eqref{tc-L2-OT} and \eqref{G-conv-Int-Y} we infer that the extension $H_{\delta,t}$ satisfies the last of the $L^2$-estimates. 
\begin{eqnarray*}
r^2\int_{X_c}i^{n^2}H_{\delta,t}\wedge_{_{h}} \overline{H}_{\delta,t}\;e^{-\varphi_{\delta,t}-\tau\psi}
&\le&
\int_{X_c}\,\frac{\;i^{n^2}H_{\delta,t}\wedge_{_{h}} \overline{H}_{\delta,t}\;e^{-\varphi_{\delta,t}-\tau\psi}\; }{\;|\sigma|^{2r}_{_H}\,(\log|\sigma|_{_H})^2\;}
\\
\\
&\le&
320\,K_r\,I_{Y_c}(f,\sigma,\varphi_{\delta,t})\,.
\end{eqnarray*}
By extracting a weak limit first in the $t$-parameter and second in the $\delta$-parameter, as permitted by claim \ref{w-lim-wei} (A), we infer the existence of a limit holomorphic extension $H_{c,\infty}$ of $f$ such that
$$
r^2\int_{X_c}i^{n^2}H_{c,\infty}\wedge_{_{h}} \overline{H}_{c,\infty}\;e^{-\varphi-\tau\psi}\le 320\,K_r\,I_{Y_c}(f,\sigma,\varphi)<+\infty\,.
$$
\\
\\
{\bf C). Construction of the global extension with uniform $L^2$-estimate.} 
\\
We set 
$\left<\cdot,\cdot\right>_{\delta,\tau,t}:=\left<\cdot,\cdot\right>_{\omega_c,h_{\delta,\tau,t}}$. According to step (B)
the section $f$ admits a holomorphic extension $H_c\equiv H_{c,\infty}$ over $X_c$ with 
\begin{eqnarray}\label{L2-Rext}
\int_{X_c}i^{n^2}H_c\wedge _h \overline{H}_c\,e^{-\varphi}<+\infty\,,
\end{eqnarray}
since $\psi$ is uniformly bounded over $X_c$.
We set $H_{c,\varepsilon}:=\theta(\varepsilon^{-2}|\sigma|^2)H_c$ and we obtain
\begin{eqnarray*}
\bar\partial_L H_{c,\varepsilon}
=
(1+\varepsilon^{-2}|\sigma|^2)\,\theta'(\varepsilon^{-2}|\sigma|^2)\,\chi'(S _{\varepsilon})^{-1}\,\bar\partial \eta_{\varepsilon}\wedge H_c\,,
\end{eqnarray*}
as in step (B).
Let $(W_{\delta,t})_{t>0}$ be a non-increasing family of relatively compact open sets in $U_{\delta}$ such that $\mu_t<\tau$ over $W_{\delta,t}$ and $U_{\delta}=\bigcup_{t>0}W_{\delta,t}$.
So by combing the lemma \ref{O-T-Key-Es} with the positivity estimate \eqref{G-Lost-Pos} as we did in \eqref{Fund-Ex-es} we infer for all $v\in L^2_{\delta,\tau,t}(X_c,\Lambda^{n,1}T^*_X\otimes L)$
\begin{eqnarray}\label{Sh-Glo-Ex-Est}
\left|\left<\bar\partial_L H_{c,\varepsilon}\,,v\right>_{\delta,\tau,t}\right|^2
\le
8\theta'(\varepsilon^{-2}|\sigma|^2)^2\,|H_c|^2_{\delta,\tau,t}\left< L^{\varepsilon}_{\delta,\tau,t}v,v\right>_{\delta,\tau,t}\,,
\end{eqnarray}
over the open set $W_{\delta,t}$, with $L^{\varepsilon}_{\delta,\tau,t}:=\left[\Theta^{\varepsilon}_{\delta,\tau,t}\,,\,\omega_c^*\right]$. 
Moreover the inequality \eqref{MRG-Lost-Pos} implies in a similar way the estimate
\begin{eqnarray}\label{Rf-Glo-Ex-Est}
\left|\left<\bar\partial_L H_{c,\varepsilon}\,,v\right>_{\delta,\tau,t}\right|^2 \le
8\theta'(\varepsilon^{-2}|\sigma|^2)^2\,|H_c|^2_{\delta,\tau,t}\left< (L^{\varepsilon}_{\delta,\tau,t}+k^2_{c,\varepsilon}I)v,v\right>_{\delta,\tau,t}\,,
\end{eqnarray}
over $X_c$.
Thus if we set 
$$
C^{\varepsilon}_{\delta,\tau,t}:=8n!\int_{X_c}\theta'(\varepsilon^{-2}|\sigma|^2)^2\,i^{n^2}H_c\wedge_{\delta,\tau,t}\overline{H}_c\,,
$$
we infer from \eqref{Sh-Glo-Ex-Est} and \eqref{Rf-Glo-Ex-Est}
\begin{eqnarray*}
\left|\left(\bar\partial_L H_{c,\varepsilon}\,,v\right)_{\delta,\tau,t}\right|^2
\le 
 2C^{\varepsilon}_{\delta,\tau,t}\left(( L^{\varepsilon}_{\delta,\tau,t}+k^2_{c,\varepsilon} Q_{\delta,t})v,v\right)_{\delta,\tau,t}\,,
\end{eqnarray*}
with $ Q_{\delta,t}$ the multiplication operator by the characteristic function of the set $X_c\smallsetminus W_{\delta,t}$. 
Let now $ P_{\delta,\tau,t}$ be the $L^2_{\delta,\tau,t}(X_c)$-orthogonal projector on the closed subspace
$$
\Ker \bar\partial_L\subset L^2_{\delta,\tau,t}(X_c,\Lambda^{n,1}T^*_X\otimes L)\,,
$$
By corollary \ref{OT-Pert-Ex} we infer the existence of a solution $(U^{\varepsilon}_{\delta,\tau,t},R^{\varepsilon}_{\delta,\tau,t})$ of the perturbed $\bar\partial$-equation 
$$
\bar\partial_LU^{\varepsilon}_{\delta,\tau,t}+k_{c,\varepsilon} P_{\delta,\tau,t}\, Q_{\delta,t}\,R^{\varepsilon}_{\delta,\tau,t} =\bar\partial_L H_{c,\varepsilon}\,,
$$
with the $L^2$-estimate
\begin{eqnarray*}
n!\int_{X_c}(\eta_{\varepsilon}+\lambda_{\varepsilon})^{-1}\,i^{n^2}U^{\varepsilon}_{\delta,\tau,t}\wedge_{\delta,\tau,t} \overline{U^{\varepsilon}_{\delta,\tau,t}}
\;+\;\int_{X_c}|R^{\varepsilon}_{\delta,\tau,t}|^2_{\delta,\tau,t}\,dV_{\omega_c}\;\le\;  2C^{\varepsilon}_{\delta,\tau,t}\,.
\end{eqnarray*}
By claims \eqref{w-lim-wei} and  \eqref{Prod-w-lim} we can extract a weak limit in the $t$-parameter so as to drop the dependence on $t$ in the last two relations, with $Q_{\delta}$ being the multiplication operator by the characteristic function of the set $V_{\delta}$. (Here we use the relation \eqref{Zer-Meas}.)
Moreover we observe
$$
\lim_{\delta\rightarrow 0} C^{\varepsilon}_{\delta,\tau}= C_{\varepsilon,\tau}:=8n!\int_{X_c}\theta'(\varepsilon^{-2}|\sigma|^2)^2\,i^{n^2}H_c\wedge_h\overline{H}_c\,|\sigma|^{-2r}\,e^{-\varphi-\tau\psi}<+\infty\,,
$$
by \eqref{L2-Rext}.
We infer by claims \eqref{w-lim-wei} (A), \eqref{Prod-w-lim} and the identity \eqref{IntesYY}, the weak convergence $Q_{\delta}\,R^{\varepsilon}_{\delta,\tau}\rightarrow 0$ as $\delta \rightarrow 0$. By claim \eqref{w-lim-wei} (B) we infer $P_{\delta,\tau} Q_{\delta}\,R^{\varepsilon}_{\delta,\tau}\rightarrow 0$ weakly as $\delta \rightarrow 0$.
By the claim \eqref{w-lim-wei} (A) we can extract a weak limit $U_{\varepsilon,\tau}$ as $\delta\rightarrow 0$, with the $L^2$-estimate
\begin{eqnarray*}
&&n!\int_{X_c}(\eta_{\varepsilon}+\lambda_{\varepsilon})^{-1}\,i^{n^2}U_{\varepsilon,\tau} \wedge_h \overline{U}_{\varepsilon,\tau}\;|\sigma|^{-2r}\,e^{-\varphi-\tau\psi}
\\
\\
&\le&
2C_{\varepsilon}:=16\,n!\int_{X_c}\theta'(\varepsilon^{-2}|\sigma|^2)^2\,i^{n^2}H_c\wedge\overline{H}_c\,|\sigma|^{-2r}\,e^{-\varphi}<+\infty\,.
\end{eqnarray*}
and solution of the $\bar\partial$-problem $\bar\partial_L U_{\varepsilon,\tau}=\bar\partial_L H_{c,\varepsilon}$. The fact that $\psi$ is uniformly bounded from above over $X_c$ allow us to extract again a weak  limit $U_{\varepsilon}$ in the $\tau$-parameter, solution of the equation $\bar\partial_L U_{\varepsilon}=\bar\partial_L H_{c,\varepsilon}$, with the $L^2$-estimate
\begin{eqnarray}
&&n!\int_{X_c}(\eta_{\varepsilon}+\lambda_{\varepsilon})^{-1}\,i^{n^2} U_{\varepsilon}\wedge_h \overline{U}_{\varepsilon}\;|\sigma|^{-2r}\,e^{-\varphi}\le  2C_{\varepsilon}
\end{eqnarray}
Moreover $U_{\varepsilon\,{|Y}}\,=\,0$ since the complex codimension of $Y$ is $r$. 
In this way we obtain a holomorphic  section $F'_{\varepsilon,c}:=H_{c,\varepsilon}-U_{\varepsilon}$ which coincides with $f$ over $Y_c$ for all $\varepsilon>0$. Then the existence of the required $L^2$-extension $F'_c$ over $X_c$ follows from the same argument explained at the end of step (B). The $L^2$-estimate allows to take again a weak limit in the $c$-parameter in order to find a $L^2$-extension $F'$ over $X_{A'}$, which in his turn extends to a global holomorphic section $F$ with the required $L^2$-estimate (with constant $C_r=320K_r$).
\\
\\
{\bf D). End of the proof.} 
We prove at this point the following elementary facts needed in step (B).
\begin{lem}\label{Conv-Int-lm}
In the setting of the theorem \ref{Ext-tm} let $F$ be an arbitrary continuous $($holomorphic$)$ extension of $f$ let $\varphi$ be a continuous $($$Y$-admissible$)$ weight. Then 
$$
\lim_{\varepsilon\rightarrow 0}\;\int_X\theta'(\varepsilon^{-2}|\sigma|^2)^2\,i^{n^2}F\wedge _h \overline{F}\;|\sigma|^{-2r}e^{-\varphi}=K_r\,I_Y(f,\sigma,\varphi)\,, 
$$
with $K_r>0$ a constant depending only on $r$ and $\theta$.
\end{lem}
$Proof$. Let
$\omega$ be an arbitrary hermitian form on $T_X$. The conclusion follows by combining the identity
$$
\frac{1}{n!}\,|f|^2_{\omega,h}\,|\Lambda^rd\sigma|^{-2}_{\omega,H}\,dV_{Y,\omega}=\frac{1}{r!}\,i^{p^2}f_{\sigma}\wedge_{_{h,H}}\bar f_{\sigma}\,,
$$
with 
\begin{eqnarray*}
\lim_{\varepsilon\rightarrow 0}\;\int_Xg\,\theta'(\varepsilon^{-2}|\sigma|^2)^2 \,|\sigma|^{-2r}e^{-\varphi}dV_{X,\omega}=
I_Y:=r!\,K_r\int_Yg\,|\Lambda^rd\sigma|^{-2}_{\omega,H}\,e^{-\varphi}dV_{Y,\omega}\,,
\end{eqnarray*}
where $g:=\frac{1}{n!}\,|F|^2_{\omega,h}$. We prove this last equality. 
Since the problem can be localised we assume that $E$ is topologically trivial over $X$. Let $\theta$ be an $H$-orthonormal trivialisation of the bundle $E$, consider $s:=\theta\circ \sigma:X\rightarrow\C^r$ and set $s_{\varepsilon}:=\varepsilon^{-1}s$,
$$
\left\{\Lambda^rd s\right\}^2_{\omega}:=\frac{n!\,r!}{p!}\;\frac{i^{r^2}\Lambda^rds\wedge \overline{\Lambda^rds}\wedge\omega^p}{\omega^n}\,.
$$ 
The fact that $\sigma$ is holomorphic implies $d\sigma=\partial\sigma$ over $Y$, which in his turn implies $ds=\partial s$ over $Y$. We infer
$$
|\Lambda^rd\sigma|^2_{\omega,H}=
|\Lambda^r\partial s|^2_{\omega}=\frac{n!\,r!}{p!}\;\frac{i^{r^2}\Lambda^r\partial s\wedge \overline{\Lambda^r\partial s}\wedge\omega^p}{\omega^n}=\left\{\Lambda^rd s\right\}^2_{\omega}\,,
$$
over $Y$.
Consider also
\begin{eqnarray*}
\alpha&:=&r!\,\theta'(\varepsilon^{-2}|z|^2)^2\,|z|^{-2r}\,i^{r^2}dz\wedge d\bar z\,,\qquad K_r:=\frac{1}{r!}\,\int_{z\in \C ^r}\alpha(z)\,,
\\
\\
\beta&:=&\frac{1}{p!}\,g\,\left\{\Lambda^rd s\right\}^{-2}_{\omega}\,\omega^p\,.
\end{eqnarray*}
We infer
\begin{eqnarray*}
J_{\varepsilon}:=\int_Xs_{\varepsilon}^*\alpha\wedge \beta\,e^{-\varphi}&=& \int_Xr!\,\theta'(\varepsilon^{-2}|s|^2)^2\,|s|^{-2r}\,i^{r^2}\Lambda^rds\wedge \overline{\Lambda^rds}\wedge\beta\,e^{-\varphi}
\\
\\
&=&
\int_Xg\,\theta'(\varepsilon^{-2}|\sigma|^2)^2 \,|\sigma|^{-2r}e^{-\varphi}dV_{X,\omega}\,.
\end{eqnarray*}
On the other hand
\begin{eqnarray*}
\lim_{\varepsilon\rightarrow 0} J_{\varepsilon}=\lim_{\varepsilon\rightarrow 0} \int_{z\in \C ^r}\alpha(z)\cdot\int_{y\in s^{-1}(\varepsilon z)}\beta(y)e^{-\varphi(y)}=
\int_{ \C ^r}\alpha\cdot\int_{ s^{-1}(0)}\beta e^{-\varphi}=I_Y\,.
\end{eqnarray*}
This is obvious in the case $\varphi$ is continuous. In the case $\varphi$ is $Y$-admissible it follows directly from the definition.\hfill$\Box$
\\
\\
We prove now the estimate \eqref{Ban-Est} needed in steps (B) and (C). 
We observe first the inequality
\begin{eqnarray*}
S_{\varepsilon}\le \log(e^{-2}+\varepsilon^2)\le -2 +O(\varepsilon^2) \,.
\end{eqnarray*}
By multiplying both sides of this inequality by $S_{\varepsilon}$ (respectively $3S_{\varepsilon}$) we obtain
\begin{eqnarray}
S_{\varepsilon}^2
&\ge& 
(-2+O(\varepsilon^2)) S_{\varepsilon}\ge (-2+O(\varepsilon^2))^2=4-O(\varepsilon^2)\label{Ban-Est-I}
\\\nonumber
\\
S_{\varepsilon}^2
&\ge& 
-6S_{\varepsilon}+O(\varepsilon^2)S_{\varepsilon}\label{Ban-Est-II}
\end{eqnarray}
Combining the identity $(\chi')^2/\chi''=(2-t)^2$ with the inequality $\eta_{\varepsilon}\le \varepsilon-2S_{\varepsilon}$ we infer
\begin{eqnarray}\label{Ban-Est-III}
\eta_{\varepsilon}+\lambda_{\varepsilon}\le 4+S_{\varepsilon}^2-6S_{\varepsilon}+\varepsilon\,.
\end{eqnarray}
We estimate the therms $4$ and $-6S_{\varepsilon}$ in \eqref{Ban-Est-III} respectively by means of the inequalities \eqref{Ban-Est-I} and \eqref{Ban-Est-II}. We obtain
$$
\eta_{\varepsilon}+\lambda_{\varepsilon}\le 5S_{\varepsilon}^2-O(\varepsilon^2)S_{\varepsilon}+O(\varepsilon)\le \left(5+O(\varepsilon^2)\right)S_{\varepsilon}^2+O(\varepsilon)\,,
$$
since $S_{\varepsilon}^2\ge -S_{\varepsilon}\ge 2-O(\varepsilon^2)$. We deduce the estimate \eqref{Ban-Est}.\hfill$\Box$

\subsection{Simplifications of the proof in the local or projective case.}\label{Loc-Proj}
We explain now how our general positivity estimate of step (A) simplifies in the local or projective case. 
In this setting we can assume that $X_c$ is a bounded pseudoconvex domain in $\C^n$. We will keep in part the notations introduced in the step (A) of the general proof.
We set 
$$
\tilde\Theta_{\delta,\tau}:=i\,{\cal C}_h(L)+i\partial\bar\partial(\varphi+\tau\psi+rS_{\delta})\,.
$$
We let $K>1$ such that $i\,{\cal C}_h(L)+i\partial\bar\partial \varphi\ge -K\omega$ and   $\tilde\Theta_{\delta,\tau}\ge -K\,\omega_c$ over $X_c$. The same type of computation done in step (A) shows
\begin{eqnarray*}
\tilde\Theta_{\delta,\tau}
&\ge&
\left(\tau-\frac{\delta^2K}{|\sigma|^2+\delta^2}\right)\omega \;+\;\frac{|\sigma|^2}{|\sigma|^2+\delta^2}\,\Theta\,,
\\
\\
\eta_{\varepsilon}\tilde\Theta_{\delta,\tau}
&\ge&
\tilde\tau_{\varepsilon,\delta}\,\omega\;+\;\frac{\;|\sigma|^2\chi'(S _{\varepsilon})\;}{|\sigma|^2+\varepsilon^2}\,\Theta
\;\ge\;
\tilde\tau_{\varepsilon,\delta}\,\omega\;+\;\chi'(S _{\varepsilon})\Phi_{\varepsilon}\,,
\end{eqnarray*}
for all $\delta\in (0,\varepsilon)$ and with 
$$
\tilde\tau_{\varepsilon,\delta}:=2\tau-\frac{\delta^2K_{\varepsilon}}{|\sigma|^2+\delta^2}\,.
$$
We define $\tilde\Theta^{\varepsilon}_{\delta,\tau,t}$ by formally replacing $\varphi$ with $\varphi_t$ in the definition of $\Theta^{\varepsilon}_{\delta,\tau}$ given in step (A) and we set
$$
\tilde V_{\delta}:=X_c\cap\left\{\tilde\tau_{\varepsilon,\delta}<\tau\right\}=X_c\cap \left\{|\sigma|^2< \delta^2(K_{\varepsilon}/\tau-1)\right\}\,.
$$
The same type of argument explained in step (A) implies that
by means of usual regularising kernels we can construct a family $(\varphi_t)_{t>0}\subset C^{\infty}(\overline{X}_c,\R)$ such that $\varphi_t\downarrow \varphi$ as $t\rightarrow 0$ over $X_c$ and
\begin{eqnarray}\label{L-Lost-Pos}
\tilde\Theta^{\varepsilon}_{\delta,\tau,t}\ge  \frac{\varepsilon^2}{2|\sigma|^2}\,i\,\partial\eta_{\varepsilon}\wedge \bar\partial\eta_{\varepsilon}+(\tau-\tilde\mu_t)\,\omega\,,
\end{eqnarray}
over the interior $\tilde U_{\delta}$ of $X_c\smallsetminus \tilde V_{\delta}$,
with $(\tilde\mu_t)_{t>0}\equiv(\tilde\mu_t^{\varepsilon,\delta,\tau})_{t>0}\subset  \R_{>0}$ a family of constants such that $\tilde\mu_t\downarrow 0$ as $t\rightarrow 0$. 
\\
\\
The rest of the proof in the local or projective case follows, with the obvious simplifications, the lines of steps (B), (C) and (D). 
\\
\\
{\bf The case $\rk_{_{\C}}E=1$.}
In the case $X$ projective (or bounded pseudoconvex domain) and $\rk_{_{\C}}E=1$, the Lelong-Poincar\'{e} formula allows to regularise with arbitrary small lost of positivity over $X_c$. More precisely by using this formula we deduce that the positivity condition $\Theta\ge 0$ is equivalent to the condition
\begin{eqnarray*}
i\,{\cal C}_h(L)+i\partial\bar\partial\varphi-i\,{\cal C}_H(E)\ge 0\,.
\end{eqnarray*}
In particular $\Theta\ge i\,{\cal C}_h(L)+i\partial\bar\partial-i\,{\cal C}_H(E)\ge 0$. By using again the Lelong-Poincar\'{e} formula we infer that the curvature condition 
$\Theta\ge \alpha^{-1}i\,{\cal C}_H(E)$ is equivalent to the condition
\begin{eqnarray*}
i\,{\cal C}_h(L)+i\partial\bar\partial\varphi-i\,{\cal C}_H(E)\ge \alpha^{-1}i\,{\cal C}_H(E)\,.
\end{eqnarray*}
So by means of usual regularising kernels we can construct a family $(\varphi_t)_{t>0}\subset C^{\infty}(\overline{X}_c,\R)$ such that $\varphi_t\downarrow \varphi$ as $t\rightarrow 0$ over $X_c$, a family of constants $(\mu'_t)_{t>0}\subset \R_{>0}$ such that $\mu'_t\downarrow 0$ as $t\rightarrow 0$ and. 
\begin{eqnarray*}
i\,{\cal C}_h(L)+i\partial\bar\partial\varphi_t-i\,{\cal C}_H(E)&\ge& -\mu'_t\,\omega\,,
\\
\\
i\,{\cal C}_h(L)+i\partial\bar\partial\varphi_t-i\,{\cal C}_H(E)&\ge& \alpha^{-1} i\,{\cal C}_H(E)-\mu'_t\,\omega\,.
\end{eqnarray*}
So if we set 
$$
\hat\Theta_{\tau,t}:=i\,{\cal C}_h(L)+i\partial\bar\partial \left(\varphi_t+\tau\psi+\log|\sigma|^2\right)\,,
$$
we infer once again by the Lelong-Poincar\'{e} formula the inequalities
\begin{eqnarray}\label{Rg1-pos-es}
\hat\Theta_{\tau,t}\ge (\tau-\mu'_t)\omega\,\qquad
\mbox{and}\qquad
\hat\Theta_{\tau,t}\ge \alpha^{-1}i\,{\cal C}_H(E)+(\tau-\mu'_t)\omega\,,
\end{eqnarray}
which in their turn imply
$$
\eta_{\varepsilon}\hat\Theta_{\tau,t}\ge \frac{|\sigma|^2\,\alpha\,\chi'(S _{\varepsilon})}{|\sigma|^2+\varepsilon^2}\hat\Theta_{\tau,t}\ge \chi'(S _{\varepsilon})\Phi_{\varepsilon}\,,
$$
for $\mu'_t<\tau$.
The fact that the extension $F_{\infty}$ constructed in step (B) of the general proof satisfies $|\bar\partial_LF_{\infty}|^2_{\omega,h}\le C|\sigma|^2$ over the set $X_c\cap\{|\sigma|<\varepsilon\}$ allows us to construct a holomorphic extension over $X_c$. In fact by applying standard $L^2$-methods one can solve the $\bar\partial$-equation $\bar\partial u_c=\bar\partial F_{\infty}$ with respect to the weight $|\sigma|^{-2r}h\,e^{-\bar\tau\psi}$ and with $\bar\tau>>0$ sufficiently big to inshur 
$$
i\,{\cal C}_h(L)+i\partial\bar\partial(\bar\tau\psi+rS_{\delta})\ge \omega\,,
$$
over $X_c$ for all $\delta\in (0,1)$. So  we obtain an extension $H'_c:=F_{\infty}-u_c$ and we set $H'_{c,\varepsilon}:=\theta(\varepsilon^{-2}|\sigma|^2)H'_c$. 
\\
\\
By using the corollary \ref{OT-Pert-Ex} with respect to the complete K\"{a}hler metric of $X_c\smallsetminus Y_c$ (see \cite{Dem3}) and with hermitian metric $|\sigma|^{-2r}h\,e^{-\varphi_t-\tau\psi}$, we infer the existence of a solution of the $\bar\partial$-problem 
$$
\bar\partial_L U^{\varepsilon}_{\tau,t}=\bar\partial_L H'_{c,\varepsilon}\,,
$$ 
over $X_c$ with the adequate $L^2$-estimate. (The solution $U^{\varepsilon}_{\tau,t}$ extends over $X_c$ thanks to a standard lemma in \cite{Dem3}). 
\\
The argument at the end of step (B) implies the existence
of a
holomorphic extension $F_{t,\tau}$ with the uniform $L^2$-estimate 
$$
J_{X_c}(F_{t,\tau},\sigma,\varphi_t+\tau\psi)\le C_r I_{Y_c}(f,\sigma,\varphi_t)\,.
$$
By claim \ref{w-lim-wei} (A) we can extract a weak limit in the $t$-parameter so as to drop it.
The fact that $\psi$ is uniformly bounded from above over $X_c$ allow us to extract again a weak  limit in the $\tau$-parameter.
\\
\\
{\bf Further developments.} 
We learn from J.-P. Demailly that his regularising technique \cite{Dem2} my admit the following more geometric version. One can expect the existence of holomorphic hermitian vector bundles $(F_k,h_k)$ and holomorphic sections $f_k\in H^0(X,F_k)$ such that the approximations in \cite{Dem2} could be written in the form
$$
\varphi_k\;=\;k^{-1}\log|f_k|^2_{h_k}+u_k\,,
$$
with $u_k$ smooth. By assuming this, one could drop the $Y$-admissibility condition on the weight $\varphi$ in the more general situation $X$ almost Stein and $\rk_{_{\C}}E=1$. In fact in this more general framework the extension result will follow by replacing the standard regularisations $\varphi_t$ with Demailly's approximations $\varphi_k$ in the positivity estimates \eqref{Rg1-pos-es}. Then the conclusion will follow by applying the corollary  \ref{OT-Pert-Ex} as just explained, to the extension $H_{c,\infty}$ constructed in step $(B)$ of the general proof with respect to the complete K\"{a}hler metric of the manifold 
$$
X_c\smallsetminus (Y_c\cup \varphi_k^{-1}(-\infty))\,.
$$
This is possible thanks to the inequality \eqref{nrm-Comp-inq}.
\\
\\
{\bf Remark.} As observed before, in the case $\rk_{_{\C}}E=1$ the curvature conditions \eqref{Cruc-Pos-hyp} are equivalent to the condition 
$$
i\,{\cal C}_h(L)+i\partial\bar\partial\varphi\;\ge\; (1+k\alpha^{-1})\,i\,{\cal C}_H(E)\,,
$$
for all $k=0,1$. In the case $\rk_{_{\C}}E=r$ arbitrary, the curvature conditions \ref{Cruc-Pos-hyp} can be replaced by a much less sharp condition.
\begin{theorem}\label{rud-Ext-tm}
Let $(L,h)$ and $(E,H)$ be two holomorphic hermitian vector bundles of rank $\rk_{_{\C}}L=1$, $\rk_{_{\C}}E=r$ over a projective variety $X$ of complex dimension $n$, let $\omega>0$ be a hermitian form over $X$, let $\lambda\equiv \lambda_{E,H}^{\omega}\in C^0(X,\R)$ such that 
$$
i{\cal C}_H(E)\le \lambda\I_{E}\otimes \omega\,,
$$
let $\sigma\in H^0(X,E)$ such that $|\sigma|_{_H}\le e^{-\alpha}$, $\alpha\in \R_{\ge 1}$ on $X$ and  set 
$$
Y:=\overline{Y}_0\,,\;\;Y_0:=\left\{x\in X\,|\,\sigma(x)=0\,,\Lambda^rd\sigma(x)\not=0\right\}\,.
$$
Consider also a quasi-plurisubharmonic function
$\varphi\in L^1(X,\R)$ such that  
the restriction $\varphi_{{|Y}}$ is  not identically $-\infty$ on any connected component of $Y$ and 
\begin{eqnarray}\label{rudim-Pos-hyp}
i\,{\cal C}_h(L)+i\partial\bar\partial\varphi\;\ge\;\lambda(r+k\alpha^{-1})\omega\,,
\end{eqnarray}
for all $k=0,1$.
Then there exist a uniform constant $C_r>0$ depending only on $r$ and on a fixed cut off function
such that for any section $f\in H^0(Y,K_X\otimes L)$ with $I_Y(f,\sigma,\varphi)<+\infty$
there exist
$F\in H^0(X,K_X\otimes L)$ such that $F=f$ over $Y$ and $J_X(F,\sigma,\varphi)\,\le \,C_rI_Y(f,\sigma,\varphi)$.
\end{theorem}
$Proof$. The curvature condition \eqref{rudim-Pos-hyp} allow to regularise $\varphi$ with arbitrary small lost of positivity over $X_c$.
Then the curvature conditions \eqref{Cruc-Pos-hyp} hold for $\hat\Theta_{\tau,t}$. This follows from the inequality 
$$
i\partial\bar\partial \log|\sigma|^2_{_H} \;\ge\; -H\left(i{\cal C}_H(E)\sigma,\sigma\right)|\sigma|^{-2}_{_H}\,.
$$
The rest of the proof follows precisely the same lines of the case $\rk_{_{\C}}E=1$ previously explained.\hfill$\Box$
\section{Singular pluri-extension results.}
In this section we combine the previous singular versions of the Ohsawa-Takegoshi-Manivel $L^2$-extension theorem with an invariance of plurigenera  technique invented by Siu \cite{Siu1}, \cite{Siu2} and drastically simplified in \cite{Pa1}.
\begin{theorem}{\bf (The codimension one case).}\label{plr-Ext-rk1-II}
Let  $(L,h_L)$ and $(E,h_E)$ be two holomorphic line bundles over a polarised projective manifold $(X,\omega)$ of complex dimension $n$. Let $\sigma\in H^0(X,E)$ and 
$\varphi_L\in L^1(X,\R)$ be a quasi-plurisubharmonic function such that 
\begin{eqnarray}\label{rk1-KLT-II}
\int_Y\frac{\;e^{-\varphi_L}\,\omega^{n-1}}{\,|d\sigma|^2_{\omega,h_E}}<+\infty\,,
\end{eqnarray}
with $Y:=\overline{Y}_0$, $Y_0:=\left\{x\in X\,|\,\sigma(x)=0\,,d\sigma(x)\not=0\right\}$
and
\begin{eqnarray}\label{rk1-Pos-hypII}
i\,{\cal C}_{h_L}(L)+i\partial\bar\partial\varphi_L\;\ge\; (1+k\alpha^{-1})\,i\,{\cal C}_H(E)\,,
\end{eqnarray}
for all $k=0,1$ and for some $\alpha\in \R_{\ge 1}$. Let also $(F,h_Fe^{-\varphi_F})$ be a pseudoeffective singular hermitian line bundle over $X$ such that  
the restriction $\varphi_{F\,{|Y}}$ is  not identically $-\infty$ on any connected component of $Y$.
Then for any $u\in H^0(Y,m(K_X+ L)+F)$ such that 
\begin{eqnarray}\label{rk1-int}
\int_Y\,\frac{|u|^2_{h_m}}{\,|d\sigma|^2_{\omega,h_E}}\,e^{-\varphi_L-\varphi_F}\,\omega^{n-1}\,<+\infty 
\end{eqnarray}
with $h_m:=\Omega_X^{-m}\otimes h^m_L\otimes h_F$ and $\Omega_X:=\omega^n/(n!)$
there exist a section
$U\in H^0(X,m(K_X+ L)+F)$ such that $U=u$ over $Y$ and
$$
\int_X\,|U|^2_{h_m}e^{-\varphi_L-\varphi_F/m}\,\omega^n\,<+\infty\,.
$$
\end{theorem}
\begin{theorem}{\bf (The arbitrary codimension case).}\label{r-Ext-tmII}
Let $(L,h)$ and $(E,h_E)$ be two holomorphic hermitian vector bundles of rank $\rk_{_{\C}}L=1$, $\rk_{_{\C}}E=r$ over a polarised projective manifold $(X,\omega)$, let $\lambda\equiv \lambda_{E,h_E}^{\omega}\in C^0(X,\R)$ such that 
$$
i{\cal C}_{h_E}(E)\le \lambda\I_{E}\otimes \omega\,.
$$  
Consider also $\sigma\in H^0(X,E)$ and 
 a quasi-plurisubharmonic function 
$\varphi_L\in L^1(X,\R)$ such that 
\begin{eqnarray}\label{rRK-KLT-II}
\int_Y\frac{\;e^{-\varphi_L}\,\omega^{n-r}}{\,|\Lambda^rd\sigma|^2_{\omega,h_E}}<+\infty\,,
\end{eqnarray}
with $Y:=\overline{Y}_0$, $Y_0:=\left\{x\in X\,|\,\sigma(x)=0\,,\Lambda^rd\sigma(x)\not=0\right\}$
and
\begin{eqnarray}\label{rudim-Pos-hyp}
i\,{\cal C}_{h_L}(L)+i\partial\bar\partial\varphi_L\;\ge\;\lambda (r+k\alpha^{-1})\omega\,,
\end{eqnarray}
for all  $k=0,1$ and for some $\alpha\in \R_{\ge 1}$. 
Let also $(F,h_Fe^{-\varphi_F})$ be a pseudoeffective singular hermitian line bundle over $X$ such that  
the restriction $\varphi_{F\,{|Y}}$ is  not identically $-\infty$ on any connected component of $Y$.
Then for any $u\in H^0(Y,m(K_X+ L)+F)$ such that
\begin{eqnarray}\label{rk-r-int}
\int_Y\,\frac{|u|^2_{h_m}}{\,|\Lambda^rd\sigma|^2_{\omega,h_E}}\,e^{-\varphi_L-\varphi_F}\,\omega^{n-r}\,<+\infty 
\end{eqnarray}
there exist
$U\in H^0(X,m(K_X+ L)+F)$ such that $U=u$ over $Y$ and
$$
\int_X  |U|^2_{h_m}e^{-\varphi_L-\varphi_F/m}\,\omega^n<+\infty\,.
$$
\end{theorem}
In the case the quasi-plurisubharmonic function $\varphi$ is with complex analytic singularities we can assume much sharp curvature conditions thanks to the main $L^2$-extension result \ref{Ext-tm}. In fact hold the following result.
\begin{theorem}{\bf(The analytic singularities case).}\label{ASG-Ext-tmII}
Let $(L,h)$ and $(E,h_E)$ be two holomorphic hermitian vector bundles of rank $\rk_{_{\C}}L=1$, $\rk_{_{\C}}E=r$ over a polarised projective manifold $(X,\omega)$. Consider also $\sigma\in H^0(X,E)$
and a quasi-plurisubharmonic function 
$\varphi_L\in L^1(X,\R)$ with complex analytic singularities such that \eqref{rRK-KLT-II} hold over $Y$ as in theorem \ref{r-Ext-tmII}.
and
\begin{eqnarray}\label{Cruc-Pos-hypII}
0\;\le\;
i\,{\cal C}_{h_L}(L)+i\partial\bar\partial \left(\varphi_L+r\log|\sigma|^2_{h_E}\right) 
\;\ge\; 
\alpha^{-1}h_E\left(i{\cal C}_{h_E}(E)\sigma,\sigma\right)|\sigma|^{-2}_{h_E}\,,
\end{eqnarray}
for some $\alpha\in \R_{\ge 1}$. Let also $(F,h_Fe^{-\varphi_F})$ be a pseudoeffective singular hermitian line bundle over $X$ such that $\varphi_F$ is with analytic singularities and
the restriction $\varphi_{F\,{|Y}}$ is  not identically $-\infty$ on any connected component of $Y$.
Then hold the conclusion of theorem \ref{r-Ext-tmII}.
\end{theorem}
{\bf Proof.} We will prove all this results at the same time. We can assume $\sup_X\varphi_L=\sup_X\varphi_F=0$.
For all $\nu\in \N$ let $k_{\nu}:=\max\{k\in \N\,:\,km\le \nu\}$, $q_{\nu}:=\nu-k_{\nu}m=0,...,m-1$  and $L_m:=m(K_X+L)+F$ equipped with hermitian metric $h_m$. We choose an ample line bundle $A$ over $X$ such that
\\
\\
{\bf (A1)} for all $q=0,...,m-1$ the line bundle $(q(K_X+L)+A)_{{|Y}}$ is base point free, globally generated by some family $(s_{q,j})_{j=1}^{N_q}\subset H^0(Y,q(K_X+L)+A)$.
\\
\\
{\bf (A2)} the restriction map $H^0(X,L_m+A)\rightarrow H^0(Y,L_m+A)$ is surjective.
\\
\\
Let $\Omega_Y$ be the volume form over $Y$ induced from the metric $\omega$.
We fix also a smooth hermitian metric $h_A$ on $A$ and
we note by $|\cdot|_{\nu}$ the norm of the induced hermitian metric 
$$
H_{\nu}:=\Omega_X^{-\nu}\otimes h_L^{\nu}\otimes h_F^{k_{\nu}}\otimes h_A
$$ 
over the line bundle ${\cal L}_{\nu}:=\nu(K_X+L)+k_{\nu}F+A$. 
The assumption (A1) implies
$$
\max_{0\le p,q\le m-1}\max_Y\;\frac{\sum_{j=1}^{N_q}|s_{q,j}|^2_q}{\sum_{t=1}^{N_p}|s_{p,t}|^2_p}\;=\;C_1\;<\;+\infty
$$
With the  notations introduced so far hold the following lemma.
\begin{lem}\label{Pau-lm}
There exist a constant $C>0$ such that for all $\nu\in \N_{\ge m}$ there exist a family of  sections 
$$
(S_{\nu,j})_{j=1}^{M_{\nu}}\subset H^0(X,{\cal L}_{\nu})\,,
$$
$M_{\nu}:=N_{q_{\nu}}$ such that $S_{\nu,j\,{|Y}}=u^{k_{\nu}}\otimes s_{q_{\nu},j}$ for all $j=1,...,M_{\nu}$ and 
$$
\int_XB_{\nu}/B_{\nu-1}\,\Omega_X\le C\,,
$$
with $B_{\nu}:=\sum_{j=1}^{M_{\nu}}|S_{\nu,j}|_{\nu}^2$ , $B_{m-1}:=1$ by convention.
\end{lem}
$Proof$. The proof goes by induction. The statement is obvious for $\nu=m$ by the assumption (A2). So we assume it true for $\nu$ and we prove it for $\nu+1$. 
Let $\delta_{\nu}=0$ if $q_{\nu}\le m-2$ and $\delta_{\nu}=1$ if $q_{\nu}=m-1$.
We have 
$$
{\cal L}_{\nu+1}=K_X+{\cal L}_{\nu}+L+\delta_{\nu}F\,,
$$ 
and $H_{\nu+1}=\Omega_X^{-1}\otimes H_{\nu}\otimes h_L\otimes \delta_{\nu}h_F$. The singular hermitian line bundle 
$$
({\cal L}_{\nu}+L+\delta_{\nu}F \,,\,H_{\nu}B_{\nu}^{-1}\otimes h_Le^{-\varphi_L}\otimes \delta_{\nu}h_Fe^{-\varphi_F})
$$ 
satisfies the analogue of the curvature conditions \eqref{rk1-Pos-hypII}, \eqref{rudim-Pos-hyp} and \eqref{Cruc-Pos-hypII} required to apply the singular versions of the $L^2$-extension theorem stated in the section \ref{sng-ext-rst}. We infer that any section $s\in H^0(Y, {\cal L}_{\nu+1})$ such that
$$
I_{Y,\,\nu+1}(s):=\int_Y\frac{|s|^2_{\nu+1}\,e^{-\varphi_L-\delta_{\nu}\varphi_F}\,\Omega_Y}{\,|\Lambda^rd\sigma|^2_{\omega,h_E}\,B_{\nu}}<+\infty\,,
$$
admits an extension 
$S\in H^0(X, {\cal L}_{\nu+1})$ 
satisfying the estimates
$$
\int_X\frac{|S|^2_{\nu+1}}{B_{\nu}}\,\Omega_X\;\le \;\int_X\frac{|S|^2_{\nu+1}}{B_{\nu}}\,e^{-\varphi_L-\delta_{\nu}\varphi_F}\,\Omega_X\;\le \;C_0\,I_{Y,\,\nu+1}(s)\,,
$$
where $C_0>0$ is a uniform constant.
We distinguish two cases.
\\
\\
{\bf Case I}. In the case $q_{\nu}\le m-2$ hold
$$
{\cal L}_{\nu+1}=k_{\nu}m(K_X+L)+(q_{\nu}+1)(K_X+L)+k_{\nu}F+A\,,
$$
This implies that the section 
$$
s:=u^{k_{\nu}}\otimes s_{q_{\nu}+1,j}\in H^0(Y,{\cal L}_{\nu+1})\,,
$$ 
$j=1,...,M_{\nu+1}$, (notice that in case I hold $q_{\nu+1}=q_{\nu}+1$). This combined with the fact that
$$
B_{\nu\,{|Y}}=\sum_{t=1}^{M_{\nu}}|u^{k_{\nu}}\otimes s_{q_{\nu},t}|_{\nu}^2\,,
$$
with the definition of the constant $C_1$ and with the assumptions \eqref{rk1-KLT-II}, \eqref{rRK-KLT-II},
allow us to apply the previous version of the Ohsawa-Takegoshi-Manivel extension theorem in order to obtain the required extensions $S_{\nu+1,j}$.
\\
\\
{\bf Case II}. In the case $q_{\nu}= m-1$ 
hold ${\cal L}_{\nu+1}=(k_{\nu}+1)L_m+A$, which implies $s:=u^{k_{\nu}+1}\otimes s_{0,j}\in H^0(Y,{\cal L}_{\nu+1})$. This combined with the fact that
$$
B_{\nu\,{|Y}}=\sum_{t=1}^{N_{m-1}}|u^{k_{\nu}}\otimes s_{m-1,t}|_{\nu}^2\,,
$$
with the definition of the constant $C_1$ and with the assumptions \eqref{rk1-int}, \eqref{rk-r-int}, allow us to apply the previous version of the Ohsawa-Takegoshi extension theorem in order to obtain the required extensions $S_{\nu+1,j}$.\hfill$\Box$
\\
\\
Let $\dashint$ be the integral mean value operator. By lemma \ref{Pau-lm} and Jensen inequality we infer
$$
\dashint_X(\log B_{\nu}-\log B_{\nu-1})\,\Omega_X\;\le\;\log\;\dashint_XB_{\nu}/B_{\nu-1}\,\Omega_X\;\le\; C'\,.
$$
Moreover the singular hermitian line bundle 
$$
\left({\cal L}_{km}\,,\,H_{km}B_{km}^{-1}\right)\equiv \left(kL_m+A\,,\, h_m^k\otimes h_A\,B_{km}^{-1}\right)\,,
$$
is pseudoeffective.
So in conclusion we got the following.
\begin{eqnarray*}
&&
\frac{1}{k}\int_X\log B_{km}\,\Omega_X \;\le\; mC'\,,
\\
\\
&&
i\,{\cal C}_{h_m}(L_m)\,+\,\frac{1}{k}\, i\,\partial\bar\partial \log B_{km} \;\ge\; -\, \frac{1}{k}\,i\,{\cal C}_{h_A}(A)\,,
\\
\\
&&
\frac{1}{k}\,\log B_{km\,{|Y}}\;=\;\log |u|^2_{h_m}\,+\,\frac{1}{k}\,\log\,\sum_{j=0}^{N_0}|s_{0,j}|^2_{h_A}\,.
\end{eqnarray*}
By well known elementary properties of quasi-plurisubharmonic functions we infer that the $L^1$-norm of the functions $\psi_k:=\frac{1}{k}\log B_{km}$ is uniformly bounded. By the $L^1$-compactness of quasi-plurisubharmonic functions we infer the existence of a subsequence $\psi_{k_j}$ convergent in the $L^1$-norm and a.e to a quasi-plurisubharmonic function $\psi$ such that 
$$
i\,{\cal C}_{h_m}(L_m)+i\,\partial\bar\partial \psi\ge 0\,,
$$ 
and $\psi_{|Y}\ge \log |u|^2_{h_m}$. This last inequality follows from the mean value inequality for plurisubharmonic functions. Thus the decomposition
$$
L_m=K_X+(m-1)(K_X+L)+L+F
$$
shows that the singular hermitian line bundle 
$$
\left((m-1)(K_X+L)+L+F \,,\, \Omega_X^{-(m-1)}\otimes h^m_L\otimes h_F \,e^{-\theta}\right)
$$
with 
$$
\theta:=\frac{m-1}{m}\,\psi+\varphi_L+\frac{1}{m}\,\varphi_F\,,
$$ 
satisfies the analogue of the curvature conditions \eqref{rk1-Pos-hypII}, \eqref{rudim-Pos-hyp} and \eqref{Cruc-Pos-hypII} required to apply the singular versions of the Ohsawa-Takegoshi-Manivel extension theorem stated in the section \ref{sng-ext-rst}.
Moreover the $L^2$ condition
\begin{eqnarray*}
&&\int_Y\frac{|u|^2_{h_m}\,e^{-\theta}}{\,|\Lambda^rd\sigma|^2_{\omega,h_E}}\,\Omega_Y
\le
\int_Y\frac{|u|^{2/m}_{h_m}\,e^{-\varphi_L-\varphi_F/m}}{\,|\Lambda^rd\sigma|^2_{\omega,h_E}}\,\Omega_Y
\\
\\
&\le&
\left(\int_Y|u|^2_{h_m}e^{-\varphi_F}\;\frac{e^{-\varphi_L}\,\Omega_Y}{\,|\Lambda^rd\sigma|^2_{\omega,h_E}}\right)^{\frac{1}{m}}
\left(\int_Y\frac{e^{-\varphi_L}\,\Omega_Y}{\,|\Lambda^rd\sigma|^2_{\omega,h_E}}\right)^{\frac{m-1}{m}}
<+\infty\,,
\end{eqnarray*}
follows from the Holder inequality applied to the finite volume measure $|\Lambda^rd\sigma|^{-2}_{\omega,h_E}e^{-\varphi_L}\,\Omega_Y$ and the $L^2$ assumptions \eqref{rk1-int}, \eqref{rk-r-int}.
Thus we apply the singular versions of the $L^2$-extension theorem of the section \ref{sng-ext-rst} in order to get the required extension $U$.\hfill$\Box$ 
\\
\\
We infer the following immediate consequences of the theorem \ref{plr-Ext-rk1-II}.
\begin{corol}\label{alg-plurI}
Let $(L,h_{_{L}}e^{-\varphi_{_{L}}})$ be a pseudoeffective line bundle over a projective manifold $X$ and let $Z\subset X$ be a hypersurface such the divisor $-Z$ is semi-ample $($i.e its stable base locus is empty$)$ and 
\begin{eqnarray*}
\int_Z\frac{\,e^{-\varphi_{_{L}}}}{\left|d\sigma_{_{Z}}\right|^2_{h_{_{Z}}}}<+\infty\,,
\end{eqnarray*}
with $\sigma_{_{Z}} \in H^0(X,{\cal O}(Z))$ such that $\div \sigma_{_{Z}}=Z$ and with $h_{_{Z}}$ an arbitrary  smooth hermitian metric on ${\cal O}(Z)$. 
Let also $(F,h_Fe^{-\varphi_F})$ be a pseudoeffective singular hermitian line bundle over $X$ such that  
the restriction $\varphi_{F\,{|Z}}$ is  not identically $-\infty$ on any connected component of $Z$.
Then for any $u\in H^0(Z,m(K_X+ L)+F)$ 
such that
$$
\int_Z\,\frac{|u|^2_{_{mL,F}}}{\,\left|d\sigma_{_{Z}}\right|^2_{h_{_{Z}}}}\,e^{-\varphi_{_{L}}-\varphi_{_{F}}}\,<+\infty\,,
$$
with $|\cdot|_{_{mL,F}}$ the norm induced by the hermitian metric $h^m_{_{L}}\otimes h_{_{F}}$,
there exist
$U\in H^0(X,m(K_X+ L)+F)$ such that $U=u$ over $Z$ and 
$$
\int_X\,|U|^2_{_{mL,F}}\,e^{-\varphi_{_{L}}-\varphi_{_{F}}/m}\,<+\infty\,.
$$
\end{corol}
\begin{corol}\label{alg-plurII}
Let $Z\subset X$ be a hypersurface inside a projective manifold $X$ and let $(L,h_{_{L}})$ be a smooth hermitian line bundle such that some positive multiple $pL$ can be decomposed as $pL=A+E$ where $A$ is a semi-ample line bundle such that $\alpha A-pZ$ is also semi-ample for some $\alpha\in \N_{\ge 1}$ and $(E,h_{_{E}}e^{-\varphi_{E}})$ is a pseudoeffective line bundle such that  
\begin{eqnarray*}
\int_Z\frac{\,e^{-\varphi_{_{E}}/p}}{\left|d\sigma_{_{Z}}\right|^2_{h_{_{Z}}}}<+\infty\,,
\end{eqnarray*}
with $\sigma_{_{Z}} \in H^0(X,{\cal O}(Z))$ such that $\div \sigma_{_{Z}}=Z$ and with $h_{_{Z}}$ an arbitrary  smooth hermitian metric on ${\cal O}(Z)$. 
Let also $(F,h_{_{F}}\,e^{-\varphi_F})$ be a pseudoeffective singular hermitian line bundle over $X$ such that  
the restriction $\varphi_{F\,{|Z}}$ is  not identically $-\infty$ on any connected component of $Z$.
Then for any $u\in H^0(Z,m(K_X+Z+ L)+F)$ such that
$$
\int_Z\,\frac{\;|u|^2_{_{m(L+Z),F}}}{\,\left|d\sigma_{_{Z}}\right|^2_{h_{_{Z}}}}\;e^{-\varphi_{_{E}}/p\,-\,\varphi_{_{F}}}\,<+\infty\,,
$$
where $|\cdot|_{_{m(L+Z),F}}$ is the norm induced by the hermitian metric $h_{_{Z}}^m\otimes h^m_{_{L}}\otimes h_{_{F}}$,
there exist
$U\in H^0(X,m(K_X+Z+L)+F)$ such that $U=u$ over $Z$ and  
$$
\int_X\,|U|^2_{_{m(L+Z),F}}\,e^{-\varphi_{_{E}}/p\,-\,\varphi_{_{F}}/m} \,<+\infty\,,
$$
In particular if $Z$ is smooth and $F$ is semi-ample then the restriction morphism
$$
H^0(X,m(K_X+Z+L)+F)\longrightarrow H^0(Z,m(K_Z+ L_{{|Z}})+F_{{|Z}})\,,
$$
is surjective.
\end{corol}
$Proof$. We choose smooth metrics $h_A$ and $h_{\alpha A-pZ}$ on $A$ and $\alpha A-pZ$ respectively with semipositive curvature. We equip ${\cal O}(Z)$ with the smooth metric induced by $h_A$ and $h_{\alpha A-pZ}$ and we apply in this setting the theorem \ref{plr-Ext-rk1-II}.\hfill$\Box$
\\
\\
This last result should be compared with the proofs and results obtained in \cite{Pa1}, \cite{Dem4}, \cite{Tak}, \cite{Ha-Mc}, \cite{Var}, \cite{Be-Pa1}, \cite{Be-Pa2}. 
\section{Appendix}
\subsection{Proof of the classic $L^2$-extension result.}
We give an essentially section \ref{sng-ext-rst} self independent prove of the extension result in the case $X\subset \C^n$ is a bounded pseudoconvex domain and $Y:=Y_0$. The hermitian vector bundles $(L,h)$ and $(E,H)$ are assumed to be trivial ($\rk_{_{\C}}E$ arbitrary) and $\varphi\in \Psh(U)$, whith $U\supset\supset X$ an open set. 
Let $(\varphi_{t})_{t>0}$ be a smooth family of plurisubharmonic functions such that $\varphi_{t}\downarrow \varphi$ as $t\rightarrow 0$. By \eqref{Hol-Pos-est} and \eqref{Part-Pos-Est} we infer
$$
\Theta_{t,\delta}^{\varepsilon}:=\eta_{\varepsilon}\,i\partial\bar\partial(\varphi_{t}+S_{\delta})-i\partial\bar\partial\eta_{\varepsilon}-\lambda^{-1}_{\varepsilon}\,i\partial \eta _{\varepsilon}\wedge \bar\partial \eta _{\varepsilon}
\,\ge \,
\frac{\varepsilon^2}{2|\sigma|^2}\,i\partial \eta _{\varepsilon}\wedge \bar\partial \eta _{\varepsilon}\,.
$$
The canonical section $f$ admits a holomorphic extension $F_{\infty}$ which can be constructed by classic $L^2$-theory methods. Consider now 
$F_{\varepsilon}:=\theta(\varepsilon^{-2}|\sigma|^2)F_{\infty}$ and the expression
\begin{eqnarray*}
\bar\partial F_{\varepsilon}=
(1+\varepsilon^{-2}|\sigma|^2)\,\theta'(\varepsilon^{-2}|\sigma|^2)\,\chi'(S _{\varepsilon})^{-1}\,\bar\partial \eta_{\varepsilon}\wedge F_{\infty}\,,
\end{eqnarray*}
Combining the lemma \ref{O-T-Key-Es} with the corollary \ref{OT-Pert-Ex} (with $\rho=0$) as we did in step (B) of the general proof, we infer the existence of a solution $U_{\varepsilon,t,\delta}$ of the $\bar\partial$-equation $\bar\partial U_{\varepsilon,t,\delta}=\bar\partial F_{\varepsilon}$, with the $L^2$-estimate
\begin{eqnarray*}
&&n!\int_X(\eta_{\varepsilon}+\lambda_{\varepsilon})^{-1}\,i^{n^2}U_{\varepsilon,t,\delta}\wedge \overline{U}_{\varepsilon,t,\delta}\,(|\sigma|^2+\delta^2)^{-r}\,e^{-\varphi_{t}}
\\
\\
&\le&
8\,n!\int_X\theta'(\varepsilon^{-2}|\sigma|^2)^2\,i^{n^2}F_{\infty}\wedge\overline{F}_{\infty}\,(|\sigma|^2+\delta^2)^{-r}\,e^{-\varphi_{t}}\,.
\end{eqnarray*}
By the claim \eqref{w-lim-wei} (A) we can extract a weak limit solution $U_{\varepsilon,t}$ as $\delta\rightarrow 0$, with the $L^2$-estimate
\begin{eqnarray*}
&&n!\int_X(\eta_{\varepsilon}+\lambda_{\varepsilon})^{-1}\,i^{n^2}U_{\varepsilon,t}\wedge \overline{U}_{\varepsilon,t}\,|\sigma|^{-2r}\,e^{-\varphi_{t}}
\\
\\
&\le&
C_{\varepsilon,t}:=8n!\int_X\theta'(\varepsilon^{-2}|\sigma|^2)^2\,i^{n^2}F_{\infty}\wedge\overline{F}_{\infty}\,|\sigma|^{-2r}\,e^{-\varphi_{t}}\,.
\end{eqnarray*}
We infer $U_{\varepsilon,t\,{|Y}}\,=\,0$ since the complex codimension of $Y$ is $r$. 
In this way we obtain a holomorphic canonical section $F_{\varepsilon,t}:=F_{\varepsilon}-U_{\varepsilon,t}$ which coincides with $f$ over $Y$ for all $\varepsilon,t>0$.
Combining \eqref{Ban-Est} with the previous inequality we infer
\begin{eqnarray}\label{Apx}
n!\int_XS_{\varepsilon}^{-2}\,i^{n^2}U_{\varepsilon,t}\wedge \overline{U}_{\varepsilon,t}\,(|\sigma|^2+\varepsilon^2)^{-r}\,e^{-\varphi_{t}}\le \left(5+O(\varepsilon)\right)C_{\varepsilon,t}\,.
\end{eqnarray}
Moreover the fact that $F_{\varepsilon}$ is bounded and has support contained in $\{|\sigma|<\varepsilon\}$ implies the estimate
\begin{eqnarray*}
n!\int_XS_{\varepsilon}^{-2}\,i^{n^2}F_{\varepsilon}\wedge \overline{F}_{\varepsilon}\,(|\sigma|^2+\varepsilon^2)^{-r}\,e^{-\varphi_{t}}
&\le&
C(\log 2\varepsilon^2)^{-2}\int_{|\sigma|<\varepsilon} \varepsilon^{-2r}dV
\\
\\
&\le&
C'(\log 2\varepsilon^2)^{-2}\,,
\end{eqnarray*}
for some uniform constants $C,C'>0$.
Combining this with \eqref{Apx} we infer the last of the following inequalities
\begin{eqnarray}
\alpha_{\mu}n!\int_X\,i^{n^2}F_{\varepsilon,t}\wedge \overline{F}_{\varepsilon,t}\,e^{-\varphi_{t}}
&\le&
\rho^2n!\int_Xi^{n^2}F_{\varepsilon,t}\wedge \overline{F}_{\varepsilon,t}\,(|\sigma|^2+\varepsilon^2)^{-(r-2\mu)}e^{-\varphi_{t}}\nonumber
\\\nonumber
\\
&\le&
n!\int_XS_{\varepsilon}^{-2}\,i^{n^2}F_{\varepsilon,t}\wedge \overline{F}_{\varepsilon,t}\,(|\sigma|^2+\varepsilon^2)^{-r}\,e^{-\varphi_{t}}\nonumber
\\\nonumber
\\
&\le&
\left(5+O(\varepsilon)\right)C_{\varepsilon,t}+C'(\log 2\varepsilon^2)^{-2}\,,\label{apx2}
\end{eqnarray}
with $\alpha_{\mu}>0$ a constant depending on $\mu\in (0,r/2)$ and $\alpha$. By lemma \ref{Conv-Int-lm} we infer the convergence
\begin{eqnarray}\label{APX-conv-Int-Y}
\lim_{\varepsilon\rightarrow 0}C_{\varepsilon,t}=n!\,8\,K_r\,I_Y(f,\sigma,\varphi_{t})<+\infty
\end{eqnarray}
By weak compactness we can extract a weakly convergent sequence $F_{\varepsilon_k,t}$ with limit $F^{t}$ as $\varepsilon_k\rightarrow 0$. 
Moreover  $F^{t}=f$ over $Y$ as follows directly from the definition of Bergman Kernel of the domain $X$.
The fact that $(-S_{\varepsilon})^{-1}$ is uniformly bounded implies by claim \ref{Prod-w-lim}
$$
(-S_{\varepsilon_k})^{-1}F_{\varepsilon_k,t}\longrightarrow (-2\log|\sigma|)^{-1}F^{t}
$$ 
weakly as $\varepsilon_k\rightarrow 0$. By claim \ref{w-lim-wei} (A), \eqref{apx2} and \eqref{APX-conv-Int-Y} we infer $F^{t}$ is the required $L^2$-extension with respect to $\varphi_{t}$. The conclusion follows by extracting a weak limit in the $t$-parameter as permitted by claim \ref{w-lim-wei} (A). Notice that in this very particular case we obtain a constant  $C_r=160K_r$ in the $L^2$-estimate. \hfill$\Box$

\vspace{1cm}
\noindent
Nefton Pali
\\
Universit\'{e} Paris Sud, D\'epartement de Math\'ematiques 
\\
B\^{a}timent 425 F91405 Orsay, France
\\
E-mail: \textit{nefton.pali@math.u-psud.fr}
\end{document}